\documentclass[12pt]{article}
\usepackage{amsmath,amsthm,amsfonts,amssymb}

\theoremstyle{plain}
\newtheorem{theorem}{Theorem}[section]
\newtheorem{corollary}{Corollary}[section]
\newtheorem{proposition}{Proposition}[section]
\newtheorem{lemma}{Lemma}[section]
\theoremstyle{definition}
\newtheorem{remark}{Remark}[section]

\begin{document}

\title{\large\bf
TRANSFORMATIONS OF ORDER ONE AND QUADRATIC FORMS
ON WIENER SPACES
}
\author{\normalsize
Setsuo TANIGUCHI}
\date{}

\maketitle

\begin{abstract}
It will be shown that transformations of order one on
the Wiener space give rise to quadratic forms as exponents
of change of variables formulas, and conversely
every exponentially integrable quadratic form has 
a transformation of order one realizing the form in such a
manner. 
Several expressions of corresponding change of variables
formulas are also discussed.
\end{abstract}

\section*{Introduction}
\label{sec:intro}

To see the purpose of this paper, we start with an
elementary observation on the $N$-dimensional Euclidean
space $\mathbb{R}^N$.
Let $A:\mathbb{R}^N\to \mathbb{R}^N$ be a linear mapping and
put $B=-(A+A^*+A^*A)$, where $A^*$ is the adjoint linear 
mapping of $A$. 
Assume that the maximal eigenvalue of $B$, say $\Lambda(B)$,
is less than one. 
Then $I+A$ is bijective (cf.~Remark~\ref{r.transf}(i)) and
it holds that 
\[
    |\det(I+A)|\int_{\mathbb{R}^N} \varphi(x+Ax)
      e^{\langle Bx,x\rangle/2} g_N(x)dx
    =\int_{\mathbb{R}^N} \varphi(x) g_N(x)dx
\] 
for every bounded and continuous 
$\varphi:\mathbb{R}^N\to \mathbb{R}$, where $I$ is the
identity mapping of $\mathbb{R}^N$,
$\langle\cdot,\cdot\rangle$ is the Euclidean inner product
in $\mathbb{R}^N$, 
$g_N(x)=(2\pi)^{-N/2}e^{-|x|^2/2}$ for $x\in \mathbb{R}^N$,
and $dx$ stands for the integration with respect to the
$N$-dimensional Lebesgue measure.
Hence the quadratic form $\langle Bx,x\rangle$,
$x\in \mathbb{R}^N$, arises from the linear transformation 
$I+A$. 
Conversely, given self-adjoint linear mapping
$B:\mathbb{R}^N\to \mathbb{R}^N$ with 
$\Lambda(B)<1$, setting 
$A=(I-B)^{1/2}-I$, where $(\cdots)^{1/2}$ is the
self-adjoint and non-negative definite
square root of linear mapping, we see that the identity 
holds, that is, the quadratic form $\langle Bx,x\rangle$,
$x\in \mathbb{R}^N$, has a transformation $I+A$ 
realizing it via the identity.
Thus linear transformations and quadratic forms corresponds
to each other bi-directionally. 
The purpose of this paper is to see that such bi-directional
relationship continues to hold on Wiener spaces, that is, 
a transformation of order one on the Wiener space, 
which is a counterpart of $I+A$ and will be defined later,
gives rise to a quadratic form on the Wiener space as the
exponent of the change of variables formula, and conversely
a quadratic form has a transformation of order one realizing
it in such a manner.   

To state our results precisely, we introduce notations.
Let $T>0$, $d\in \mathbb{N}$, and $\mathcal{W}$ be the space
of continuous functions from $[0,T]$ to $\mathbb{R}^d$
vanishing at $0$. 
Denote by $\mathcal{L}_2$ the space of
$\mathbb{R}^{d\times d}$-valued square integrable
functions on $[0,T]^2$ with respect to the two-dimensional
Lebesgue measure, where $\mathbb{R}^{d\times d}$ is the
space of $d\times d$ real matrices.
The norm in $\mathcal{L}_2$ is given by
\[
    \|\kappa\|_2
    =\biggl(\int_0^T \int_0^T |\kappa(t,s)|^2 ds dt
     \biggr)^{1/2}
    \quad\text{for }\kappa\in\mathcal{L}_2,
\]
where $|\cdot|$ is the Euclidean norm in 
$\mathbb{R}^{d\times d}$.
We write as $\kappa=\kappa^\prime$ in $\mathcal{L}_2$ 
for $\kappa,\kappa^\prime\in \mathcal{L}_2$ if 
$\|\kappa-\kappa^\prime\|_2=0$.
Denote by $\mathcal{H}$ the Cameron-Martin subspace of
$\mathcal{W}$, that is, the space of absolutely continuous
$h\in \mathcal{W}$ possessing a square integrable derivative 
$h^\prime=((h^\prime)^1,\dots,(h^\prime)^d)$.
$\mathcal{H}$ is a real separable Hilbert space with the
inner product  
\[
    \langle h,g\rangle_{\mathcal{H}}
    =\int_0^T \langle h^\prime(t),g^\prime(t) \rangle dt
    \quad\text{for }h,g\in \mathcal{H},
\]
where $\langle\cdot,\cdot\rangle$ is the Euclidean inner
product in $\mathbb{R}^d$.
For $\kappa=(\kappa_j^i)_{1\le i,j\le d}\in\mathcal{L}_2$,
define the Wiener functional 
$F_\kappa:\mathcal{W}\to \mathcal{H}$ by
\begin{equation}\label{eq:f.kappa}
    \langle F_\kappa,h\rangle_{\mathcal{H}}
    =\sum_{i,j=1}^d \int_0^T \biggl(\int_0^T \kappa_j^i(t,s)
     d\theta^j(s)\biggr)(h^\prime)^i(t)dt
    =\int_0^T \biggl\langle 
         \int_0^T \kappa(t,s)d\theta(s),h^\prime(t)
         \biggr\rangle dt
\end{equation}
for $h\in \mathcal{H}$, where 
(i)~$\{\theta(t)=(\theta^1(t),\dots,\theta^d(t))
 \}_{t\in[0,T]}$ is the coordinate process on $\mathcal{W}$,
that is, $\theta(t)(w)=w(t)$ for $t\in[0,T]$ and 
$w\in \mathcal{W}$, 
(ii) each $d\theta^j(t)$ is the It\^o integral with respect
to $\{\theta^j(t)\}_{t\in[0,T]}$ and 
$d\theta(t)=(d\theta^1(t),\dots,d\theta^d(t))$, and 
(iii)~in the last term we have used the matrix notation that  
each  element of $\mathbb{R}^d$ is thought of as a column
vector and $\mathbb{R}^{d\times d}$ acts on $\mathbb{R}^d$
from left, i.e.,  
$Mx=\bigl(\sum_{j=1}^d M_j^ix^j\bigr)_{1\le i\le d}$ for
$M=(M_j^i)_{1\le i,j\le d}\in \mathbb{R}^{d\times d}$ 
and $x=(x^1,\dots,x^d)\in \mathbb{R}^d$.
Since each component of 
$\int_0^T \kappa(t,s) d\theta(s)$ belongs to the Wiener
chaos of order one, letting $\iota$ be the identity mapping 
of $\mathcal{W}$, we call $\iota+F_\kappa$ 
{\it a transformation of order one}.
Denote by $\mathcal{S}_2$ the space of all
$\eta\in\mathcal{L}_2$ with the property that
$\eta(t,s)^\dagger=\eta(s,t)$ for  
$(t,s)\in[0,T]^2$, where $M^\dagger$ is the transpose of
$M\in \mathbb{R}^{d\times d}$.
For $\eta=(\eta_j^i)_{1\le i,j\le d}\in\mathcal{S}_2$,
define the Wiener functional 
$\mathfrak{q}_\eta:\mathcal{W}\to \mathbb{R}$ by
\begin{equation}\label{eq:q.eta}
    \mathfrak{q}_\eta
    =\sum_{i,j=1}^d \int_0^T \biggl(
          \int_0^t \eta_j^i(t,s)d\theta^j(s)\biggr)
          d\theta^i(t)
    =\int_0^T \biggl\langle 
          \int_0^t \eta(t,s)d\theta(s),d\theta(t)
          \biggr\rangle.
\end{equation}
Since every element of the Wiener chaos of order two is of
the form $\mathfrak{q}_\eta$ for some $\eta\in\mathcal{S}_2$
(\cite{nualart}), we call $\mathfrak{q}_\eta$ 
{\it a quadratic form}.  
For $\kappa\in\mathcal{L}_2$, define the Hilbert-Schmidt
operator $B_\kappa:\mathcal{H}\to \mathcal{H}$ and
$\eta(\kappa)\in\mathcal{S}_2$ by
\begin{align}
    & \langle B_\kappa h,g\rangle_{\mathcal{H}}
    =\int_0^T \biggl\langle \int_0^T \kappa(t,s)
      h^\prime(s)ds,g^\prime(t)\biggr\rangle dt
    \quad\text{for } h,g\in \mathcal{H},
 \label{eq:b.kappa}
    \\
    & \eta(\kappa)(t,s)
      =-\biggl\{\kappa(t,s)+\kappa(s,t)^\dagger
     +\int_0^T \kappa(u,t)^\dagger\kappa(u,s)du
     \biggr\}
    \quad\text{for } (t,s)\in[0,T]^2.
 \label{eq:eta.kappa}
\end{align}
Note that $B_\kappa=B_{\kappa^\prime}$ if 
$\kappa=\kappa^\prime$ in $\mathcal{L}_2$.
Put
\[
    \Lambda(B)
    =\sup_{\|h\|_{\mathcal{H}}=1}\langle Bh,h
         \rangle_{\mathcal{H}}
\]
for self-adjoint Hilbert-Schmidt operator
$B:\mathcal{H}\to \mathcal{H}$, where
$\|h\|_{\mathcal{H}} 
 =\sqrt{\langle h,h\rangle_{\mathcal{H}}}$.

The first aim of this paper is to show the following change
of variables formula associated with $\iota+F_\kappa$.
If $\Lambda(B_{\eta(\kappa)})<1$, then it holds that 
\begin{equation}\label{eq:transf}
    |{\det}_2(I+B_\kappa)|
    \int_{\mathcal{W}} f(\iota+F_\kappa)
    e^{\mathfrak{q}_{\eta(\kappa)}} d\mu
    =e^{\|\kappa\|_2^2/2} \int_{\mathcal{W}} f d\mu
\end{equation}
for every $f\in C_b(\mathcal{W})$,
where ${\det}_2$ stands for the regularized determinant
(\cite{GGK}), $I$ is the identity mapping of $\mathcal{H}$
to itself, $\mu$ is the Wiener measure on $\mathcal{W}$, and  
$C_b(\mathcal{W})$ is the space of
bounded and continuous $\mathbb{R}$-valued functions on
$\mathcal{W}$. 
See Theorem~\ref{t.transf}.
Thus the quadratic form $\mathfrak{q}_{\eta(\kappa)}$ arises
from the transformation of order one $\iota+F_\kappa$.
It should be noted that $e^{\mathfrak{q}_\eta}$ is
integrable with respect to $\mu$ if and only if 
$\Lambda(B_\eta)<1$ (Lemma~\ref{l.q.eta.int}(ii)) and
hence the assumption that $\Lambda(B_{\eta(\kappa)})<1$ is
the best possible one for the identity \eqref{eq:transf} to 
hold.
Further we shall show that the transformation
$\iota+F_\kappa$ has an inverse transformation
$\iota+F_{\widehat{\kappa}}$ with 
$\widehat{\kappa}\in \mathcal{L}_2$ so that
\eqref{eq:transf} turns into the identity
\begin{equation}\label{eq:inv.transf}
    |{\det}_2(I+B_\kappa)| \int_{\mathcal{W}} f
    e^{\mathfrak{q}_{\eta(\kappa)}} d\mu
    =e^{\|\kappa\|_2^2/2} \int_{\mathcal{W}} 
     f(\iota+F_{\widehat{\kappa}}) d\mu
\end{equation}
for every $f\in C_b(\mathcal{W})$.
See Theorem~\ref{t.inv.transf}.
To construct the inverse transformation, the functional
analytic mechanism with the help of Malliavin calculus plays
a key role. 

The second aim of this paper is to see the converse
assertion that each $\eta\in\mathcal{S}_2$ with the property
that $\Lambda(B_\eta)<1$ has a
$\kappa\in\mathcal{L}_2$ such that 
$\eta(\kappa)=\eta$ in $\mathcal{L}_2$.
In particular, 
$\mathfrak{q}_\eta=\mathfrak{q}_{\eta(\kappa)}$,
because 
$\int_{\mathcal{W}}
 |\mathfrak{q}_\eta-\mathfrak{q}_{\eta^\prime}|^2d\mu
 =\|\eta-\eta^\prime\|_2^2$ 
for $\eta,\eta^\prime\in \mathcal{S}_2$. 
Thus the quadratic form $\mathfrak{q}_\eta$ admits a
transformation of order one $\iota+F_\kappa$ which realizes
it via \eqref{eq:transf} with $\mathfrak{q}_\eta$ for
$\mathfrak{q}_{\eta(\kappa)}$.
We shall give a method to find such a $\kappa$, and evaluate
Wiener integrals 
$\int_{\mathcal{W}} f e^{\mathfrak{q}_\eta} d\mu$ of Laplace
transformation type with the help of \eqref{eq:inv.transf}.
Such a $\kappa$ is obtained by solving the equation
$(I+B_\kappa)^2=I-B_\eta$.
See Theorem~\ref{t.surj.1st}.
This equation has a unique solution in a class of
Hilbert-Schmidt operators, and hence the mapping
$\kappa\mapsto\eta(\kappa)$ is bijective if we restrict the 
domain (Proposition~\ref{p.inj.1st}).
As will be seen in Remark~\ref{r.transf}(ii), the
restriction is indispensable for the injectivity.
As another derivative of the identity \eqref{eq:transf}, we
also investigate linear transformations studied
by Cameron and Martin \cite{CM}.
Then the quadratic form $\mathfrak{q}_{\eta(\kappa)}$ has a
more explicit expression (Theorem~\ref{t.transf.cm}).

The correspondence between transformations of order one and
quadratic forms was pointed out by the author in 
\cite{st-kjm-2024}.
The transformation considered there has a 
$\kappa\in\mathcal{L}_2$ with $\kappa(t,s)=0$ for 
$0\le t<s\le T$.
The appearance of $\eta(\kappa)$ for general 
$\kappa\in \mathcal{L}_2$ and the bijectivity of the mapping
$\kappa\mapsto\eta(\kappa)$ are achieved newly in this
paper.  
Thus this paper is a goal of the investigation of the
bi-directional correspondence between transformations of
order one and quadratic forms. 

The organization of this paper is as follows.
In Section~\ref{sec:transf}, we shall see that a quadratic
form arises from a transformation of order one by proving
the change of variables formula \eqref{eq:transf}.
In Section~\ref{sec:inv.transf}, it will be seen that every
transformation $\iota+F_\kappa$ with
$\Lambda(B_{\eta(\kappa)})<1$ has an inverse transformation
of order one, and the identity \eqref{eq:inv.transf} holds.
Section~\ref{sec:surj} is devoted to showing the existence
of $\kappa$ with $\eta(\kappa)=\eta$ in $\mathcal{L}_2$ for
given $\eta$. 
In the section, evaluations of Laplace transformations of 
quadratic forms and the bijectivity of the mapping
$\kappa\mapsto\eta(\kappa)$ are also
investigated.
Furthermore the evaluation will be applied to Wiener
functionals generalizing the square norm of sample path of
the one-dimensional Wiener process. 
In Section~\ref{sec:lin.transf}, linear transformations are 
studied and a more explicit expression for
$\mathfrak{q}_{\eta(\kappa)}$ will be shown.
A comparison with Cameron and Martin's result will be also
discussed.

\section{Transformation of order one}
\label{sec:transf}

The aim of this section is to see that a quadratic form
arises from a transformation of order one as follows. 

\begin{theorem}\label{t.transf}
Let $\kappa\in\mathcal{L}_2$ and assume that
$\Lambda(B_{\eta(\kappa)})<1$.
Then \eqref{eq:transf} holds, i.e.,
\[
    |{\det}_2(I+B_\kappa)| \int_{\mathcal{W}} f(\iota+F_\kappa)
    e^{\mathfrak{q}_{\eta(\kappa)}} d\mu
    =e^{\|\kappa\|_2^2/2} \int_{\mathcal{W}} f d\mu
\]
for every $f\in C_b(\mathcal{W})$.
\end{theorem}

\begin{remark}\label{r.transf}
(i) It holds that
\begin{equation}\label{r.transf.1}
    B_{\eta(\kappa)}
    =-(B_\kappa+B_\kappa^*+B_\kappa^*B_\kappa),
\end{equation}
where $B_\kappa^*$ is the adjoint operator of $B_\kappa$.
This implies that
\[
    \inf_{\|h\|_{\mathcal{H}}=1} 
     \|(I+B_\kappa)h\|_{\mathcal{H}}^2
    =\inf_{\|h\|_{\mathcal{H}}=1} \bigl\{
       1-\langle B_{\eta(\kappa)}h,h
           \rangle_{\mathcal{H}} \bigr\}
    =1-\Lambda(
         B_{\eta(\kappa)}). 
\]
Hence the condition that 
$\Lambda(B_{\eta(\kappa)})<1$ is equivalent to that
\[
    \inf_{\|h\|_{\mathcal{H}}=1} 
     \|(I+B_\kappa)h\|_{\mathcal{H}}>0.
\]
Thus, if $\|B_\kappa\|_{\text{\rm op}}<1$, where 
$\|\cdot\|_{\text{\rm op}}$ stands for the operator norm,
then 
$\Lambda(B_{\eta(\kappa)})<1$.
\\
(ii)
Take $b,c\in \mathbb{R}$ with $b^2+c^2>0$ and 
orthonormal $h_1,h_2\in \mathcal{H}$.
Define $\kappa_1,\kappa_2\in\mathcal{L}_2$ by 
\begin{align*}
    & \kappa_1(t,s)
      =b\{h_1^\prime(s)\otimes h_1^\prime(t)
        +h_2^\prime(s)\otimes h_2^\prime(t)\}, 
    \\
    & \kappa_2(t,s)
      =c\{h_1^\prime(s)\otimes h_2^\prime(t)
       -h_2^\prime(s)\otimes h_1^\prime(t)\}
    \quad\text{for }(t,s)\in[0,T]^2, 
\end{align*}
where
\[
    x\otimes y=(y^ix^j)_{1\le i,j\le d}
    =\begin{pmatrix}
       y^1x^1 & \dots & y^1x^d \\
       \vdots & & \vdots \\
       y^dx^1 & \dots & y^d x^d
     \end{pmatrix}
\]
for 
$x=(x^1,\dots,x^d),y=(y^1,\dots,y^d)\in \mathbb{R}^d$.
By a direct computation, we see that 
$\|\kappa_1-\kappa_2\|_2^2=2(b^2+c^2)>0$, and 
\begin{align*}
   &  \eta(\kappa_1)(t,s)
    =-(2b+b^2)
       \{h_1^\prime(s)\otimes h_1^\prime(t)
        +h_2^\prime(s)\otimes h_2^\prime(t)\}, 
    \\
   & \eta(\kappa_2)(t,s)
    =-c^2\{h_1^\prime(s)\otimes h_1^\prime(t)
        +h_2^\prime(s)\otimes h_2^\prime(t)\}
    \quad\text{for }(t,s)\in[0,T]^2.
\end{align*}
Suppose that $1+c^2=(1+b)^2$.
Then $\eta(\kappa_1)=\eta(\kappa_2)$.
Further 
\[
    \langle B_{\eta(\kappa_i)}h,h\rangle_{\mathcal{H}}
    =-c^2\{\langle h_1,h\rangle_{\mathcal{H}}^2
           +\langle h_2,h\rangle_{\mathcal{H}}^2\} 
    \quad\text{for $h\in \mathcal{H}$ and $i=1,2$},
\]
and hence 
$\Lambda(B_{\eta(\kappa_i)})\le0<1$ for $i=1,2$.
Thus the mapping 
$\mathcal{L}_2\ni\kappa
 \mapsto\eta(\kappa)\in\mathcal{S}_2$
is not injective, even if the assumption in the theorem is
fulfilled.
\\
(iii) 
By Lemma~\ref{l.q.eta.int} below,
$e^{\mathfrak{q}_{\eta(\kappa)}}$ is integrable with respect
to $\mu$ if and only if $\Lambda(B_{\eta(\kappa)})<1$.
\end{remark}

To prove the theorem, we prepare lemmas.
For a real separable Hilbert space $E$, 
let $L^p(\mu;E)$ be the space of $p$th integrable $E$-valued
Wiener functionals with respect to $\mu$.
$L^p(\mu;\mathbb{R})$ is denoted as $L^p(\mu)$ simply.
Let $\mathbb{D}^\infty(E)$ be the space of infinitely
$\mathcal{H}$-differentiable $E$-valued Wiener functionals
in the sense of Malliavin calculus, whose
$\mathcal{H}$-derivatives of all orders are $p$th integrable
with respect to $\mu$ for every $p\in(1,\infty)$. 
The $\mathcal{H}$-derivative and its adjoint are written as
$D$ and $D^*$, respectively.
Both 
$D:\mathbb{D}^\infty(E)\to
   \mathbb{D}^\infty(\mathcal{H}\otimes E)$ and
$D^*:\mathbb{D}^\infty(\mathcal{H}\otimes E)\to
   \mathbb{D}^\infty(E)$ are continuous,
where $\mathcal{H}\otimes E$ is the Hilbert space of
Hilbert-Schmidt operators from $\mathcal{H}$ to $E$.
For details, see \cite{mt-cambridge}.
Put 
\[
    \mathcal{H}^{\otimes 2}=\mathcal{H}\otimes \mathcal{H}
    \quad\text{and}\quad
    \mathcal{S}(\mathcal{H}^{\otimes 2})
     =\{A\in\mathcal{H}^{\otimes 2};A^*=A\}.
\]

\begin{lemma}\label{l.b.kappa}
The mapping 
$\mathcal{L}_2\ni\kappa
 \mapsto B_\kappa\in\mathcal{H}^{\otimes 2}$ 
is bijective, where the injectivity means that 
$\kappa_1=\kappa_2$ in $\mathcal{L}_2$ if 
$B_{\kappa_1}=B_{\kappa_2}$.
If $\eta\in\mathcal{S}_2$, then 
$B_\eta\in \mathcal{S}(\mathcal{H}^{\otimes 2})$, and 
the restricted mapping 
$\mathcal{S}_2\ni\eta\mapsto 
 B_\eta\in \mathcal{S}(\mathcal{H}^{\otimes 2})$ 
is also bijective.
\end{lemma}

\begin{proof}
Since 
$\|B_{\kappa_1}-B_{\kappa_2}\|_{\mathcal{H}^{\otimes 2}}
 =\|\kappa_1-\kappa_2\|_2$ 
for $\kappa_1,\kappa_2\in\mathcal{L}_2$,
where $\|\cdot\|_{\mathcal{H}^{\otimes 2}}$ is the
Hilbert-Schmidt norm in $\mathcal{H}^{\otimes 2}$, the
mapping $\kappa\mapsto B_\kappa$ is injective. 
Let $B\in\mathcal{H}^{\otimes 2}$.
Take an orthonormal basis 
$\{h_n\}_{n=1}^\infty$ of $\mathcal{H}$.
For $N\in \mathbb{N}$, define $\kappa_N\in\mathcal{L}_2$ by 
\[
    \kappa_N(t,s)=\sum_{n,m=1}^N 
     \langle Bh_n,h_m \rangle_{\mathcal{H}} 
     [h_n^\prime(s)\otimes h_m^\prime(t)]
    \quad\text{for } (t,s)\in[0,T]^2,
\]
where $x\otimes y\in \mathbb{R}^{d\times d}$ for 
$x,y\in \mathbb{R}^d$ is defined as  
in Remark~\ref{r.transf}(ii).
Since
\[
    \|\kappa_N-\kappa_M\|_2^2
    =\sum_{M<n\vee m\le N} 
      \langle Bh_n,h_m\rangle_{\mathcal{H}}^2
    \quad\text{for }1\le M<N,
\]
where $n\vee m=\max\{n,m\}$, and 
\[
    \sum_{n,m=1}^\infty 
      \langle Bh_n,h_m\rangle_{\mathcal{H}}^2
     =\|B\|_{\mathcal{H}^{\otimes 2}}^2<\infty, 
\]
there is a $\kappa\in \mathcal{L}_2$ with 
$\|\kappa_N-\kappa\|_2\to 0$ as $N\to\infty$. 
Hence 
$\|B_{\kappa_N}-B_\kappa\|_{\mathcal{H}^{\otimes 2}}
 =\|\kappa_N-\kappa\|_2\to0$ as $N\to\infty$.
Further, if we denote by $\pi_N$ the orthogonal projection 
of $\mathcal{H}$ onto the subspace spanned by
$h_1,\dots,h_N$, then $B_{\kappa_N}=\pi_N B\pi_N$ and 
hence $\|B_{\kappa_N}-B\|_{\mathcal{H}^{\otimes 2}}\to0$ as
$N\to\infty$.
Thus  $B=B_\kappa$.
Therefore the surjectivity of the mapping follows.

For $\kappa\in\mathcal{L}_2$, define 
$\kappa^*\in\mathcal{L}_2$ by 
$\kappa^*(t,s)=\kappa(s,t)^\dagger$
for $(t,s)\in[0,T]^2$.
Since $B_\kappa^*=B_{\kappa^*}$, 
$B_\eta\in \mathcal{S}(\mathcal{H}^{\otimes 2})$ 
if $\eta\in\mathcal{S}_2$.
The restricted mapping $\eta\mapsto B_\eta$ inherits the
injectivity from the original one.
Let $B\in \mathcal{S}(\mathcal{H}^{\otimes 2})$.
By the surjectivity of the original mapping, there is an
$\kappa\in\mathcal{L}_2$ with $B_\kappa=B$.
Then $B_{\kappa^*}=B_\kappa^*=B$.
Defining $\eta\in \mathcal{S}_2$ by 
$\eta=(1/2)\{\kappa+\kappa^*\}$, we have that $B_\eta=B$.
Thus the restricted mapping is surjective.
\end{proof}

\begin{lemma}\label{l.q.eta}
Let $\eta\in\mathcal{S}_2$.
Then it holds that
\[
    \mathfrak{q}_\eta=\frac12 (D^*)^2B_\eta, 
\]
where the right term is defined by regarding 
$B_\eta\in \mathcal{S}(\mathcal{H}^{\otimes 2})$ as a 
constant function in 
$\mathbb{D}^\infty(\mathcal{H}^{\otimes 2})$ and applying
$D^*$ twice. 
Further, for any orthonormal basis $\{h_n\}_{n=1}^\infty$ of
$\mathcal{H}$, it holds that
\[
    \mathfrak{q}_\eta=\frac12\sum_{n,m=1}^\infty 
     \langle B_\eta h_n,h_m\rangle_{\mathcal{H}}
     \{(D^*h_n)(D^*h_m)-\delta_{nm}\},
\]
where the series converges in $L^p(\mu)$
for every $p\in(1,\infty)$ and $D^*h_n$ is defined
by thinking of $h_n$ as a constant function in
$\mathbb{D}^\infty(\mathcal{H})$ and applying $D^*$.
\end{lemma}

\begin{proof}
The first identity was shown in
\cite[Lemma~2.2]{st-kjm-2024}.
To see the second identity, develop 
$B_\eta\in\mathcal{H}^{\otimes 2}$ as  
\[
    B_\eta=\sum_{n,m=1}^\infty \langle B_\eta h_n,h_m
       \rangle_{\mathcal{H}} h_n\otimes h_m,
\]
where 
$h_n\otimes  h_m\in\mathcal{H}^{\otimes 2}$ is defined by
\[
    (h_n\otimes  h_m)h
    =\langle h_n,h\rangle_{\mathcal{H}} h_m
    \quad\text{for }h\in \mathcal{H}.
\]
Put
\[
    B_\eta^{(N)}=\sum_{n,m=1}^N \langle B_\eta h_n,h_m
         \rangle_{\mathcal{H}} h_n\otimes h_m
    \quad\text{for }N\in \mathbb{N}.
\]
Then $B_\eta^{(N)}\to B_\eta$ in 
$\mathcal{H}^{\otimes 2}$ as $N\to\infty$.
Due to the continuity of $D^*$,
$(D^*)^2B_\eta^{(N)}\to (D^*)^2B_\eta$ in $L^p(\mu)$ for
every $p\in(1,\infty)$.
Recall that
$(D^*)^2(h_n\otimes h_m)=(D^*h_n)(D^*h_m)-\delta_{nm}$
(\cite[(5.7.2)]{mt-cambridge}).
Thus we have that
\[
    (D^*)^2B_\eta^{(N)}=\sum_{n,m=1}^N 
      \langle B_\eta h_n,h_m\rangle_{\mathcal{H}} 
      \{(D^*h_n)(D^*h_m)-\delta_{nm}\}.
\]
In conjunction with the first identity for
$\mathfrak{q}_\eta$, this implies the second identity and
the convergence of the series in any $L^p(\mu)$.
\end{proof}

\begin{lemma}\label{l.f.kappa}
Let $\kappa\in\mathcal{L}_2$.
Then it holds that
\[
    F_\kappa=D^*B_\kappa.
\]
Further, for any orthonormal basis $\{h_n\}_{n=1}^\infty$ of
$\mathcal{H}$, it holds that
\[
    F_\kappa=\sum_{n,m=1}^\infty 
     \langle B_\eta h_n,h_m\rangle_{\mathcal{H}}
     (D^*h_n)h_m,
\]
where the series converges in
$L^p(\mu;\mathcal{H})$ for every $p\in(1,\infty)$. 
\end{lemma}

\begin{proof}
It follows from \eqref{eq:f.kappa} that
\[
    \langle D F_\kappa,h\otimes g
          \rangle_{\mathcal{H}^{\otimes 2}}
    =\langle D\langle F_\kappa,g\rangle_{\mathcal{H}},
       h\rangle_{\mathcal{H}}
    =\int_0^T \biggl\langle 
         \int_0^T \kappa(t,s)h^\prime(s)ds,g^\prime(t)
         \biggr\rangle dt
    =\langle B_\kappa h,g\rangle_{\mathcal{H}}
\]
for $h,g\in \mathcal{H}$, where 
$\langle\cdot,\rangle_{\mathcal{H}^{\otimes 2}}$ is the
Hilbert-Schmidt inner product in $\mathcal{H}^{\otimes 2}$.
Thus 
\begin{equation}\label{l.f.kappa.21}
   DF_\kappa=B_\kappa.
\end{equation}

For $h,g\in \mathcal{H}$ and 
$\Phi\in \mathbb{D}^\infty(\mathbb{R})$, it holds that 
\begin{align*}
   & \int_{\mathcal{W}} \langle D(D^*B_\kappa),
            h\otimes g\rangle_{\mathcal{H}^{\otimes 2}}
          \Phi d\mu 
     =\int_{\mathcal{W}} \langle B_\kappa,
           D\bigl((D^*(\Phi h))\otimes g\bigr)
            \rangle_{\mathcal{H}^{\otimes 2}} d\mu
   \\
   & \quad
     =\int_{\mathcal{W}} \langle 
          B_\kappa\bigl(D(D^*(\Phi h))\bigr),
           g\rangle_{\mathcal{H}}
          d\mu
     =\int_{\mathcal{W}} \langle \Phi h,
          D(D^*(B_\kappa^*g))\rangle_{\mathcal{H}}
          d\mu
   \\
   & \quad
     =\int_{\mathcal{W}} \langle h,
          B_\kappa^*g \rangle_{\mathcal{H}}\Phi 
          d\mu
     =\int_{\mathcal{W}} \langle B_\kappa,h\otimes g
          \rangle_{\mathcal{H}^{\otimes 2}} \Phi d\mu,
\end{align*}
where, to see the fourth equality, we have used the fact
that $D(D^*h)=h$ for $h\in \mathcal{H}$
(\cite[(5.1.9)]{mt-cambridge}).
Hence 
\begin{equation}\label{l.f.kappa.22}
    D(D^*B_\kappa)=B_\kappa.
\end{equation}

By \eqref{l.f.kappa.21} and \eqref{l.f.kappa.22},
we have that
$D(F_\kappa-D^*B_\kappa)=0$.
Applying \cite[Proposition~5.2.9]{mt-cambridge}, we obtain
that 
\[
    F_\kappa-D^*B_\kappa
    =\int_{\mathcal{W}}(F_\kappa-D^*B_\kappa)d\mu.
\]
Due to \eqref{eq:f.kappa}, 
$\int_{\mathcal{W}}F_\kappa d\mu=0$.
If we think of $h\in \mathcal{H}$ as a constant function
in $\mathbb{D}^\infty(\mathcal{H})$, then 
$D h=0$.
Hence it holds that
\[
    \biggl\langle \int_{\mathcal{W}} (D^* B_\kappa)d\mu,
      h\biggr\rangle_{\mathcal{H}} 
    =\int_{\mathcal{W}} \langle D^* B_\kappa,
            h\rangle_{\mathcal{H}} d\mu
    =\int_{\mathcal{W}} \langle B_\kappa,
               D h\rangle_{\mathcal{H}^{\otimes 2}}
       d\mu 
    =0
    \quad\text{for } h\in \mathcal{H}.
\]
Thus $\int_{\mathcal{W}}(D^*B_\kappa)d\mu=0$.
Therefore $F_\kappa-D^*B_\kappa=0$, that is, the first
identity holds.

To see the second identity, put
\[
    B_\kappa^{(N)}=\sum_{n,m=1}^N \langle B_\kappa h_n,
         h_m\rangle_{\mathcal{H}} h_n\otimes h_m
    \quad\text{for }N\in \mathbb{N}.
\]
Then $B_\kappa^{(N)}\to B_\kappa$ in 
$\mathcal{H}^{\otimes 2}$ as $N\to\infty$.
By the continuity of $D^*$, 
$D^*B_\kappa^{(N)}\to D^*B_\kappa$ in
$L^p(\mu;\mathcal{H})$ for every $p\in(1,\infty)$.
Since
$(D^*)(h_n\otimes h_m)=(D^*h_n)h_m$
(\cite[(5.7.3)]{mt-cambridge}), 
\[
    D^*B_\kappa^{(N)}=\sum_{n,m=1}^N \langle B_\kappa h_n,
     h_m\rangle_{\mathcal{H}} (D^*h_n)h_m.
\]
Thus we obtain the second identity and the convergence of
the series in any $L^p(\mu;\mathcal{H})$.
\end{proof}

\begin{lemma}\label{l.q.eta.int}
Let $\eta\in\mathcal{S}_2$.
\\
{\rm(i)}
Assume that $\Lambda(B_\eta)<1$.
Then it holds that
\[
    \int_{\mathcal{W}} e^{\mathfrak{q}_\eta} d\mu
    \le \exp\biggl(\frac12\biggl\{\frac12
          +\frac{0\vee\Lambda(B_\eta)}{
           3\{1-(0\vee\Lambda(B_\eta))\}^3}
         \biggr\}\|\eta\|_2^2 \biggr).
\]
{\rm(ii)}
$e^{\mathfrak{q}_\eta}\in L^1(\mu)$ if and only if 
$\Lambda(B_\eta)<1$.
\end{lemma}

\begin{proof}
Let $\eta\in \mathcal{S}_2$.
\\
(i)
By Lemma~\ref{l.b.kappa}, 
$B_\eta\in \mathcal{S}(\mathcal{H}^{\otimes 2})$.
Hence there is an orthonormal basis $\{h_n\}_{n=1}^\infty$
of $\mathcal{H}$ such that 
$B_\eta=\sum_{n=1}^\infty a_n h_n\otimes h_n$, where 
$a_n\in \mathbb{N}$ for $n\in \mathbb{N}$ and 
$\sum_{n=1}^\infty a_n^2<\infty$.
By Lemma~\ref{l.q.eta}, there is an increasing sequence
$\{N_m\}_{m=1}^\infty\subset \mathbb{N}$ such that
\[
    \mu\biggl(\frac12 \sum_{n=1}^{N_m}a_n\{(D^*h_n)^2-1\}
    \to \mathfrak{q}_\eta
    \text{ as $m\to\infty$}\biggr)=1.
\]
Since $\sup_{n\in \mathbb{N}}a_n=\Lambda(B_\eta)<1$ and 
$\{D^*h_n\}_{n=1}^\infty$ is an i.i.d.~sequence of
random variables obeying the standard normal distribution
$N(0,1)$, by Fatou's lemma, we have that
\begin{align}
    \int_{\mathcal{W}} e^{\mathfrak{q}_\eta} d\mu
    & \le \liminf_{m\to\infty} \int_{\mathcal{W}}
         \exp\biggl(\frac12\sum_{n=1}^{N_m} 
           a_n\{(D^*h_n)^2-1\}\biggr) d\mu
 \nonumber
    \\
    &   = \liminf_{m\to\infty} \prod_{n=1}^{N_m}
          \biggl(\int_{\mathbb{R}} 
            e^{(a_n/2)\{x^2-1\}} 
            \frac1{\sqrt{2\pi}} e^{-x^2/2} dx
          \biggr)
 \nonumber
    \\
    &  =\biggl\{\limsup_{m\to\infty}
          \prod_{n=1}^{N_m} (1-a_n)e^{a_n}
            \biggr\}^{-1/2}.
 \label{l.q.eta.int.21}
\end{align}

The Taylor expansion of $\log(1-a)$ about $a=0$ implies that 
\[
    \log(1-a)+a+\frac{a^2}2
    =-\int_0^a\int_0^b\int_0^c \frac2{(1-u)^3} du dc db
    \ge -\frac{0\vee a}{3\{1-(0\vee a)\}^3}\,a^2
    \quad\text{for }a<1.
\]
Since $a_n\le\Lambda(B_\eta)<1$ for 
$n\in \mathbb{N}$, we obtain that  
\begin{align*}
    \prod_{n=1}^{N_m} (1-a_n)e^{a_n}
    & =\exp\biggl(\sum_{n=1}^{N_m} \biggl\{
        \log(1-a_n)+a_n+\frac{a_n^2}2\biggr\}\biggr)
       \exp\biggl(-\frac12\sum_{n=1}^{N_m} a_n^2\biggr)
    \\
    & \ge \exp\biggl(
          -\frac{0\vee\Lambda(B_\eta)}{
            3\{1-(0\vee\Lambda(B_\eta))\}^3}
         \sum_{n=1}^{N_m} a_n^2\biggr)
       \exp\biggl(-\frac12\sum_{n=1}^{N_m} a_n^2\biggr).
\end{align*}
Thus we have that
\[
    \limsup_{m\to\infty} 
     \prod_{n=1}^{N_m} (1-a_n)e^{a_n}
     \ge \exp\biggl(-\biggl\{\frac12
       +\frac{0\vee\Lambda(B_\eta)}{
          3\{1-(0\vee\Lambda(B_\eta))\}^3}
        \biggr\} \|B_\eta\|_{\mathcal{H}^{\otimes 2}}^2\biggr).
\]
Since $\|B_\eta\|_{\mathcal{H}^{\otimes 2}}=\|\eta\|_2$,
plugging this lower estimation into \eqref{l.q.eta.int.21},
we obtain the desired upper estimation. 
\\
(ii)
By the above observation, it suffices to show that 
\begin{equation}\label{l.q.eta.int.22}
    \int_{\mathcal{W}} e^{\mathfrak{q}_\eta} d\mu
    =\infty
    \quad\text{if }\Lambda(B_\eta)\ge1.
\end{equation}
To see this, we continue to use the same development of 
$B_\eta$ as above.
Assume that $\Lambda(B_\eta)\ge1$.
Then 
$\sup_{n\in \mathbb{N}} a_n=\Lambda(B_\eta)\ge1$.
Since $\lim_{n\to\infty}a_n=0$, there is an 
$n_0\in \mathbb{N}$ such that $a_{n_0}\ge1$.
Define $\eta^\prime\in \mathcal{S}_2$ so that
$B_{\eta^\prime}=\sum_{n\ne n_0} a_n h_n\otimes h_n$.
By Lemma~\ref{l.q.eta}, we see that
\[
    \mathfrak{q}_{\eta^\prime}
    =\frac12 \sum_{n\ne n_0} a_n\{(D^*h_n)^2-1\}
    \quad\text{and}\quad
    \mathfrak{q}_\eta
    =\frac12 a_{n_0}\{(D^*h_{n_0})^2-1\}
     +\mathfrak{q}_{\eta^\prime}.
\]
The first identity yields that $D^*h_{n_0}$ and
$\mathfrak{q}_{\eta^\prime}$ are independent.
Since $D^*h_{n_0}$ obeys the standard normal distribution
$N(0,1)$, we have that 
\[
    \int_{\mathcal{W}} e^{\mathfrak{q}_\eta} d\mu
    =\biggl(\int_{\mathbb{R}}
       e^{(a_{n_0}/2)\{x^2-1\}} \frac1{\sqrt{2\pi}}
       e^{-x^2/2} dx\biggr)
      \int_{\mathcal{W}} e^{\mathfrak{q}_{\eta^\prime}} 
       d\mu
    =\infty.
\]
Thus \eqref{l.q.eta.int.22} holds.
\end{proof}

Under the identification of $\mathcal{H}$ and its dual space 
$\mathcal{H}^*$, we think of the dual space
$\mathcal{W}^*$ of $\mathcal{W}$ as a subspace of
$\mathcal{H}$.
In particular, the inclusions 
$\mathcal{W}^*\subset \mathcal{H}^*=\mathcal{H}
 \subset \mathcal{W}$ are continuous and dense.

\begin{proof}[Proof of Theorem~\ref{t.transf}]
Take an orthonormal basis $\{\ell_n\}_{n=1}^\infty$ of
$\mathcal{H}$ such that $\ell_n\in \mathcal{W}^*$ for every 
$n\in \mathbb{N}$.
For $N\in \mathbb{N}$, denote by $\pi_N$ the orthogonal 
projection of $\mathcal{H}$ onto the subspace spanned by
$\ell_1,\dots,\ell_N$.
Define $\Pi^{(N)}, \kappa^{(N)}\in\mathcal{L}_2$ by
\begin{align*}
    & \Pi^{(N)}(t,s)=\sum_{n=1}^N 
        \ell_n^\prime(s)\otimes \ell_n^\prime(t),
    \\
    & \kappa^{(N)}(t,s)=\int_0^T \int_0^T
       \Pi^{(N)}(t,u) \kappa(u,v)\Pi^{(N)}(v,s)
       dudv
    \quad\text{for $(t,s)\in[0,T]^2$.}
\end{align*}
Then $B_{\kappa^{(N)}}=\pi_NB_\kappa \pi_N$.
By this and \eqref{r.transf.1}, it holds that
\[
    \|B_{\kappa^{(N)}}-B_\kappa
        \|_{\mathcal{H}^{\otimes 2}}\to0
    \quad\text{and}\quad
    \|B_{\eta(\kappa^{(N)})}-B_{\eta(\kappa)}
         \|_{\mathcal{H}^{\otimes 2}}\to0
    \quad\text{as }N\to\infty. 
\]
Due to Lemmas~\ref{l.q.eta} and \ref{l.f.kappa} and the
continuity of $D^*$, we know that
\[
    \mathfrak{q}_{\eta(\kappa^{(N)})}\to
     \mathfrak{q}_{\eta(\kappa)}
    ~~\text{in }L^p(\mu)
    \quad\text{and}\quad
    F_{\kappa^{(N)}}\to F_\kappa
    ~~\text{in }L^p(\mu;\mathcal{H}) 
    \quad\text{as }N\to\infty
\]
for any $p\in(1,\infty)$.
Further, since
\[
    |\Lambda(B)-\Lambda(B^\prime)|
    \le \|B-B^\prime\|_{\mathcal{H}^{\otimes 2}}
    \quad\text{for }
    B,B^\prime\in \mathcal{S}(\mathcal{H}^{\otimes 2}),
\]
there is an $N_0\in\mathbb{N}$ such that
\begin{equation}\label{t.transf.21}
    \sup_{N\ge N_0}
       \Lambda(B_{\eta(\kappa^{(N)})})
    <1.
\end{equation}
Take a $p\in(1,\infty)$ such that
\[
    \alpha_p\equiv  \sup_{N\ge N_0}
     \Lambda(
         B_{p \eta(\kappa^{(N)})})
     =p \sup_{N\ge N_0} 
       \Lambda(
          B_{\eta(\kappa^{(N)})})<1.
\]
By Lemma~\ref{l.q.eta.int}, we know that
\begin{align*}
    & \sup_{N\ge N_0} \int_{\mathcal{W}} 
        \exp(p\mathfrak{q}_{\eta(\kappa^{(N)})})
        d\mu 
      =\sup_{N\ge N_0} \int_{\mathcal{W}} 
         \exp(\mathfrak{q}_{p\eta(\kappa^{(N)})})
         d\mu 
    \\
    & \quad
      \le \exp\biggl(\frac12\biggl\{\frac12
      +\frac{0\vee \alpha_p}{
         3\{1-(0\vee \alpha_p)\}^3}
      \biggr\} p^2
      \sup_{N\ge N_0}\|\eta(\kappa^{(N)})\|_2^2\biggr)
    <\infty.
\end{align*}
Hence the family 
$\{\exp(\mathfrak{q}_{\eta(\kappa^{(N)})});
 N\ge N_0\}$ 
is uniform integrable. 
Thus, on account of the continuity of the mapping 
$\mathcal{H}^{\otimes 2}\ni A\mapsto{\det}_2(I+A)$
(\cite[Theorem~XI.2.2]{GGK}), \eqref{eq:transf} follows once
we have shown the identity 
\begin{equation}\label{t.transf.22}
   |{\det}_2(I+B_{\kappa^{(N)}})| \int_{\mathcal{W}} 
        f(\iota+F_{\kappa^{(N)}})  
        \exp(\mathfrak{q}_{\eta(\kappa^{(N)})})
        d\mu 
      =e^{\|\kappa^{(N)}\|_2^2/2} \int_{\mathcal{W}} f d\mu
\end{equation}
for every $f\in C_b(\mathcal{W})$ and $N\ge N_0$.

Let $N\ge N_0$ and set 
\[
    A^{(N)}=\bigl(\langle B_\kappa \ell_n,\ell_m
      \rangle_{\mathcal{H}}\bigr)_{1\le n,m\le N}
      \in \mathbb{R}^{N\times N}.
\]
Since $B_{\kappa^{(N)}}=\pi_NB_\kappa \pi_N$, 
$B_{\kappa^{(N)}}=\sum_{n,m=1}^N
 A_{nm}^{(N)}\ell_n\otimes \ell_m$,
where $A^{(N)}=\bigl(A_{nm}^{(N)}\bigr)_{1\le n,m\le N}$.
By \eqref{r.transf.1}, we have that 
\begin{equation}\label{t.transf.23}
    A^{(N)}+A^{(N)\dagger}+A^{(N)}A^{(N)\dagger}
     =\Bigl(-\langle B_{\eta(\kappa^{(N)})}\ell_n,
         \ell_m\rangle_{\mathcal{H}}\Bigr)_{1\le n,m\le N}.
\end{equation}
If we set $\varepsilon_N=1-\Lambda(B_{\kappa^{(N)}})$, then
by \eqref{t.transf.21}, $\varepsilon_N>0$ and 
$|(I+A^{(N)\dagger})x|^2\ge \varepsilon_N|x|^2$ for every
$x\in \mathbb{R}^N$, where $I$ is the $N$-dimensional
identity matrix.
Thus the linear mapping 
$\mathbb{R}^N\ni x\mapsto 
 (I+A^{(N)\dagger})x\in \mathbb{R}^N$ is bijective.
Applying the change of variables formula for this mapping,
we obtain that
\begin{align}
    & |\det(I+A^{(N)\dagger})| 
      \int_{\mathbb{R}^N} \varphi(x+A^{(N)\dagger}x)
      \frac1{\sqrt{2\pi}^N}
      e^{-|x+A^{(N)\dagger}x|^2/2} dx
 \nonumber
    \\
    & \quad
    =\int_{\mathbb{R}^N} \varphi(x) 
      \frac1{\sqrt{2\pi}^N} e^{-|x|^2/2} dx
\label{t.transf.24}
\end{align}
for every $\varphi\in C_b(\mathbb{R}^N)$.

By \eqref{t.transf.23} and \eqref{r.transf.1}, we see that 
\begin{align*}
    & |x+A^{(N)\dagger}x|^2-|x|^2
    \\
    & \quad
      = -\sum_{n,m=1}^N 
        \langle B_{\eta(\kappa^{(N)})}\ell_n,
            \ell_m\rangle_{\mathcal{H}} 
        \{x^nx^m-\delta_{nm}\}
      +2\sum_{n=1}^N \langle B_{\kappa^{(N)}}\ell_n,
       \ell_n\rangle_{\mathcal{H}}
      +\|B_{\kappa^{(N)}}\|_{\mathcal{H}^{\otimes 2}}^2
\end{align*}
for $x=(x^1,\dots,x^N)\in \mathbb{R}^N$.
Recalling that ${\det}_2(I+M)=\det(I+M)e^{-\text{\rm tr} M}$
for $M\in \mathbb{R}^{N\times N}$, we have that
\begin{align}
    &\det(I+A^{(N)\dagger})e^{-|x+A^{(N)\dagger}x|^2/2}
     \nonumber
    \\
    & \quad
    ={\det}_2(I+A^{(N)}) 
     \exp\biggl(\frac12 \sum_{n,m=1}^N
       \langle B_{\eta(\kappa^{(N)})}\ell_n,
         \ell_m\rangle_{\mathcal{H}}
      \{x^nx^m-\delta_{nm}\}\biggr)
      e^{-|x|^2/2}e^{-\|\kappa^{(N)}\|_2^2/2}
 \label{t.transf.25}
\end{align}
for every $x=(x^1,\dots,x^N)\in \mathbb{R}^N$.

Since $B_{\kappa^{(N)}}=\pi_NB_\kappa \pi_N$ and 
$B_{\eta(\kappa^{(N)})}
 =\pi_NB_{\eta(\kappa^{(N)})} \pi_N$, 
by Lemmas~\ref{l.q.eta} and \ref{l.f.kappa}, we have
that 
\begin{align*}
    & F_{\kappa^{(N)}}=\sum_{n,m=1}^N 
       \langle B_\kappa \ell_n,\ell_m\rangle_{\mathcal{H}}
       (D^*\ell_n)\ell_m,
    \\
    & \mathfrak{q}_{\eta(\kappa^{(N)})}
      =\frac12 \sum_{n,m=1}^N 
       \langle B_{\eta(\kappa^{(N)})} \ell_n,
              \ell_m\rangle_{\mathcal{H}}
       \{(D^*\ell_n)(D^*\ell_m)-\delta_{nm}\}.
\end{align*}
Put 
$\boldsymbol{\ell}^{(N)}=(\ell_1,\dots,\ell_N):
  \mathcal{W}\to \mathbb{R}^N$ and
$D^*\boldsymbol{\ell}^{(N)}
 =(D^*\ell_1,\dots,D^*\ell_N)$.
By the above expression of $F_{\kappa^{(N)}}$, 
\[
    \ell_m(F_{\kappa^{(N)}})
    =\langle \ell_m,F_{\kappa^{(N)}}\rangle_{\mathcal{H}}
    =\sum_{n=1}^N\langle  B_\kappa\ell_n,
    \ell_m\rangle_{\mathcal{H}}(D^*\ell_n)
    \quad\text{for }1\le m\le N.
\]
Hence 
\[
    \boldsymbol{\ell}^{(N)}(F_{\kappa^{(N)}})
    =A^{(N)\dagger}D^*\boldsymbol{\ell}^{(N)}.
\]
Since $D^*\boldsymbol{\ell}^{(N)}$ obeys the standard
normal distribution $N(0,I)$ on $\mathbb{R}^N$, 
due to \eqref{t.transf.25}, 
the above expression of $\mathfrak{q}_{\eta(\kappa^{(N)})}$, 
and the fact that $D^*\ell_n=\ell_n$ for
$1\le n\le N$ (\cite[(5.1.5)]{mt-cambridge}), the
identity \eqref{t.transf.24} is rewritten as 
\[
    |{\det}_2(I+B_{\kappa^{(N)}})|\int_{\mathcal{W}} 
     [\varphi\circ \boldsymbol{\ell}^{(N)}]
     (\iota+F_{\kappa^{(N)}}) 
     \exp(\mathfrak{q}_{\eta(\kappa^{(N)})})
     d\mu
    =e^{\|\kappa^{(N)}\|_2^2/2}
     \int_{\mathcal{W}} 
     [\varphi\circ \boldsymbol{\ell}^{(N)}]
     d\mu.
\]
Thus \eqref{t.transf.22} holds for 
$f=\varphi\circ \boldsymbol{\ell}^{(N)}$.
Since both $F_{\kappa^{(N)}}$ and
$\mathfrak{q}_{\eta(\kappa^{(N)})}$ depend only on
$D^*\boldsymbol{\ell}^{(N)}$, due to the splitting property
of the Wiener measure, \eqref{t.transf.22} holds for every  
$f\in C_b(\mathcal{W})$. 
\end{proof}

\section{Inverse transformation of order one}
\label{sec:inv.transf}
The aim of this section is to show that every transformation
of order one $\iota+F_\kappa$ for $\kappa\in \mathcal{L}_2$
with $\Lambda(B_{\eta(\kappa)})<1$ has an inverse
transformation. 
If $\kappa\in \mathcal{L}_2$ and 
${\det}_2(I+B_\kappa)\ne0$, then $I+B_\kappa$ has a
continuous inverse (\cite[Theorem\,XII.1.1]{GGK}).
Rewriting $(I+B_\kappa)^{-1}-I$ as 
$-(I+B_\kappa)^{-1}B_\kappa$,
we see that 
$(I+B_\kappa)^{-1}-I\in \mathcal{H}^{\otimes 2}$
(\cite[Theorem\,IV.7.1]{GGK}).
Applying Lemma~\ref{l.b.kappa}, we find a 
$\widehat{\kappa}\in \mathcal{L}_2$ with 
\[
    B_{\widehat{\kappa}}=(I+B_\kappa)^{-1}-I.
\]
For $\kappa\in \mathcal{L}_2$ with
$\Lambda(B_{\eta(\kappa)})<1$, it follows from
\eqref{eq:transf} for $f=1$ that ${\det}_2(I+B_\kappa)\ne0$. 
Hence $\widehat{\kappa}$ is defined for such a $\kappa$.
Since
$\int_{\mathcal{W}}
  \|F_\kappa-F_{\kappa^\prime}\|_{\mathcal{H}}^2 d\mu
 =\|\kappa-\kappa^\prime\|_2^2$, 
$F_\kappa=F_{\kappa^\prime}$ if 
$\kappa=\kappa^\prime$ in $\mathcal{L}_2$.
By the injectivity given in Lemma~\ref{l.b.kappa}, 
$F_{\widehat{\kappa}}$ is determined by $\kappa$ uniquely up
to null sets.
The transformation $\iota+F_{\widehat{\kappa}}$ is an
inverse transformation of $\iota+F_\kappa$ as follows. 

\begin{theorem}\label{t.inv.transf}
Let $\kappa\in \mathcal{L}_2$ and assume that
$\Lambda(B_{\eta(\kappa)})<1$.
Then 
\begin{equation}\label{t.inv.transf.1}
    (\iota+F_\kappa)\circ(\iota+F_{\widehat{\kappa}})
    =(\iota+F_{\widehat{\kappa}})\circ(\iota+F_\kappa)
    =\iota
\end{equation}
and \eqref{eq:inv.transf} holds, i.e.,
\[
    |{\det}_2(I+B_\kappa)| \int_{\mathcal{W}} 
      f e^{\mathfrak{q}_{\eta(\kappa)}} d\mu
    =e^{\|\kappa\|_2^2/2} \int_{\mathcal{W}} 
      f(\iota+F_{\widehat{\kappa}}) d\mu
\]
for every $f\in C_b(\mathcal{W})$.
\end{theorem}

\begin{remark}\label{r.inv.transf}
(i) 
Precisely speaking, \eqref{t.inv.transf.1} means that
\[
    w+F_{\widehat{\kappa}}(w)
      +F_\kappa\bigl(w+F_{\widehat{\kappa}}(w)\bigr)
   =w+F_\kappa(w)
      +F_{\widehat{\kappa}}\bigl(w+F_\kappa(w)\bigr)
   =w
\]
for $\mu$-almost every $w\in \mathcal{W}$.
Since both $F_\kappa$ and $F_{\widehat{\kappa}}$ are
determined up to null sets, to determine terms 
$F_\kappa(w+F_{\widehat{\kappa}}(w))$
and $F_{\widehat{\kappa}}(w+F_\kappa(w))$,  we need nice
modifications of $F_\kappa$ and $F_{\widehat{\kappa}}$.
Such nice modifications will be obtained in
Lemma~\ref{l.h.inv} below.
\\
(ii)
It may be interesting to prove functional analytically
that ${\det}_2(I+B_\kappa)\ne 0$ if
$\Lambda(B_{\eta(\kappa)})<1$.
To see this, observe that 
\[
    \inf_{\|h\|_{\mathcal{H}}=1}
      \|(I-B_{\eta(\kappa)}) h\|_{\mathcal{H}}
    \ge \inf_{\|h\|_{\mathcal{H}}=1}
      \langle (I-B_{\eta(\kappa)})h,h
         \rangle_{\mathcal{H}}
    \ge 1-\Lambda(B_{\eta(\kappa)}).
\]
This implies that $I-B_{\eta(\kappa)}$ has a continuous
inverse. 
In fact, the above lower estimation yields the inequality
\[
    \|(I-B_{\eta(\kappa)}) h\|_{\mathcal{H}}\ge
    \{1-\Lambda(B_{\eta(\kappa)})\}\|h\|_{\mathcal{H}}
    \quad\text{for }h\in \mathcal{H}.
\]
This inequality implies the injectivity of
$I-B_{\eta(\kappa)}$ and  the closedness of the range of
$I-B_{\eta(\kappa)}$.  
If $h\in \mathcal{H}$ is perpendicular to the range, then,
by the self-adjointness of $I-B_{\eta(\kappa)}$ and its
injectivity, $h=0$.  
Thus $I-B_{\eta(\kappa)}$ is bijective.
Thanks to the inverse mapping theorem, it has a continuous
inverse. 
By \eqref{r.transf.1}, we have that
\[ 
    {\det}_2\bigl((I+B_\kappa^*)(I+B_\kappa)\bigr)
    ={\det}_2(I-B_{\eta(\kappa)})\ne0
\]
(\cite[Theorem\,XII.1.1]{GGK}).
Observe that
\begin{equation}\label{r.inv.transf.1}
    {\det}_2\bigl((I+B^*)(I+B)\bigr)
    ={\det}_2(I+B){\det}_2(I+B^*)
     e^{-\text{\rm tr}(B^*B)}
    \quad\text{for }B\in \mathcal{H}^{\otimes 2}.
\end{equation}
In fact, if $B$ is of trace class, then this follows from
the identities
\begin{align*}
    & {\det}_2(I+A_1)=\det(I+A_1)e^{-\text{\rm tr}A_1},
    \\
    & \det\bigl((I+A_1)(I+A_2)\bigr)
      =\det(I+A_1)\det(I+A_2)
    \quad\text{for $A_1,A_2\in \mathcal{H}^{\otimes 2}$
               of trace class}
\end{align*}
(\cite[Theorem~IX.2.1 and IV.(5.10)]{GGK}).
Then it is extended to $\mathcal{H}^{\otimes 2}$ by a
standard approximation argument.
Thus we obtain that ${\det}_2(I+B_\kappa)\ne0$.
\end{remark}

For the proof, we prepare lemmas.
A measurable set $X\subset \mathcal{W}$ is said to be
$\mathcal{H}$-invariant if $X+h\equiv\{w+h;w\in X\}$
coincides with $X$ for every $h\in \mathcal{H}$.
As is easily seen, 
$X$ is $\mathcal{H}$-invariant if and only if 
$X+h\subset X$ for every $h\in \mathcal{H}$.

\begin{lemma}\label{l.h.inv}
Let $\kappa\in\mathcal{L}_2$.
Then there is an $\mathcal{H}$-invariant set $X_\kappa$ such
that $\mu(X_\kappa)=1$ and 
\begin{equation}\label{l.h.inv.1}
    F_\kappa(w+h)=F_\kappa(w)+B_\kappa h
    \quad\text{for every $w\in X_\kappa$ and 
          $h\in \mathcal{H}$}.
\end{equation}
\end{lemma}

This lemma asserts the existence of both
$\mathcal{H}$-invariant set and nice modification of
$F_\kappa$.

\begin{proof}
Let $\{h_n\}_{n=1}^\infty$ be an orthonormal basis of
$\mathcal{H}$.
There are $\mathcal{H}$-invariant sets $X^{(n)}$, 
$n\in \mathbb{N}$, such that
\[
    \mu(X^{(n)})=1
    \quad\text{and}\quad
    (D^*h_n)(w+h)=(D^*h_n)(w)
     +\langle h_n,h\rangle_{\mathcal{H}}
\]
for every $w\in X^{(n)}$ and $h\in\mathcal{H}$ for each
$n\in \mathbb{N}$ (\cite[Lemma~5.7.7]{mt-cambridge}).
Put
\[
    X_\kappa=\biggl\{w\in \bigcap_{n=1}^\infty X^{(n)};
     \lim_{N_1,N_2\to\infty} \biggl\|
      \sum_{N_1<n\vee m\le N_2} 
      \langle B_\kappa h_n, h_m\rangle_{\mathcal{H}}
      (D^*h_n)(w)h_m\biggr\|_{\mathcal{H}}=0
    \biggr\}.
\]
For $N\in \mathbb{N}$ and $h\in \mathcal{H}$, 
if $\pi_N$ stands for the orthogonal projection of
$\mathcal{H}$ onto the subspace spanned by $h_1,\dots,h_N$,  
then
\[
    \sum_{n,m=1}^N 
     \langle B_\kappa h_n,h_m\rangle_{\mathcal{H}}
     \langle h,h_n\rangle_{\mathcal{H}} h_m
     =(\pi_N B_\kappa \pi_N) h, 
\]
and hence the sum converges to $B_\kappa h$ in $\mathcal{H}$
as $N\to\infty$.  
Thus $X_\kappa$ is $\mathcal{H}$-invariant.
By Lemma~\ref{l.f.kappa}, $\mu(X_\kappa)=1$.
On account of the same lemma, if we define the modification
of $F_\kappa$ by
\[
    F_\kappa(w)
    =\begin{cases}
       \displaystyle
       \lim_{N\to\infty} \sum_{n,m=1}^N
          \langle B_\kappa h_n, h_m\rangle_{\mathcal{H}}
          (D^*h_n)(w)h_m
          & \text{if }w\in X_\kappa,
       \\
       0 & \text{otherwise},
     \end{cases}
\]
then it satisfies \eqref{l.h.inv.1}.  
\end{proof}

\begin{lemma}\label{l.d*a}
Let $A\in\mathcal{H}^{\otimes 2}$ and 
$g\in \mathcal{H}$.
Then it holds that
$\langle D^*A,g\rangle_{\mathcal{H}}=D^*(A^*g)$.
\end{lemma}

\begin{proof}
For every $\Phi\in \mathbb{D}^\infty(\mathbb{R})$, we have
that 
\begin{align*}
    \int_{\mathcal{W}} 
       \langle D^*A,g\rangle_{\mathcal{H}}
       \Phi d\mu
    & =\int_{\mathcal{W}} 
       \langle A,D\Phi\otimes g
          \rangle_{\mathcal{H}^{\otimes 2}} d\mu
      =\int_{\mathcal{W}} 
       \langle A(D\Phi),g\rangle_{\mathcal{H}} d\mu
    \\
    & =\int_{\mathcal{W}} 
       \langle D\Phi,A^*g\rangle_{\mathcal{H}} d\mu
      =\int_{\mathcal{W}} \Phi(D^*(A^*g)) d\mu.
\end{align*}
This implies the desired identity.
\end{proof}

\begin{proof}[Proof of Theorem~\ref{t.inv.transf}]
By Lemma~\ref{l.h.inv} for $\kappa$ and $\widehat{\kappa}$,
we have that 
\begin{align*}
    & \bigl((\iota+F_\kappa)\circ
        (\iota+F_{\widehat{\kappa}})\bigr)(w)
      =w+F_\kappa(w)+F_{\widehat{\kappa}}(w)
       +B_\kappa F_{\widehat{\kappa}}(w),
    \\
    & \bigl((\iota+F_{\widehat{\kappa}})
        \circ(\iota+F_\kappa)\bigr)(w)
      =w+F_\kappa(w)+F_{\widehat{\kappa}}(w)
       +B_{\widehat{\kappa}}F_\kappa(w)
\end{align*}
for every $w\in X_\kappa\cap X_{\widehat{\kappa}}$.
By Lemmas~\ref{l.f.kappa} and \ref{l.d*a} and the identity
\[
    (I+B_\kappa)(I+B_{\widehat{\kappa}})=I
     =(I+B_{\widehat{\kappa}})(I+B_\kappa), 
\]
we see that 
\begin{align*}
    \langle
     F_\kappa+F_{\widehat{\kappa}}+B_\kappa F_{\widehat{\kappa}},
       g\rangle_{\mathcal{H}}
    & =\langle D^*B_\kappa,g\rangle_{\mathcal{H}}
       +\langle D^*B_{\widehat{\kappa}},
         g\rangle_{\mathcal{H}}
       +\langle D^*B_{\widehat{\kappa}},
         B_\kappa^* g\rangle_{\mathcal{H}}
    \\
    & =D^*\Bigl(B_\kappa^*g+B_{\widehat{\kappa}}^*g
        +B_{\widehat{\kappa}}^*(B_\kappa^*g)\Bigr)
    \\
    & =D^*\bigl[\bigl(
        \{(I+B_\kappa)(I+B_{\widehat{\kappa}})\}^*
         -I\bigr)g\bigr]
      =0
\end{align*}
and 
\[
    \langle
       F_\kappa+F_{\widehat{\kappa}}+B_{\widehat{\kappa}}F_\kappa,
       g\rangle_{\mathcal{H}}
      =D^*\bigl[\bigl(
        \{(I+B_{\widehat{\kappa}})(I+B_\kappa)\}^*
         -I\bigr)g\bigr]
      =0
\]
for any $g\in \mathcal{H}$.
Thus \eqref{t.inv.transf.1} holds.

To show \eqref{eq:inv.transf}, notice that \eqref{eq:transf} 
continues to hold for any bounded and measurable 
$\phi:\mathcal{W}\to \mathbb{R}$, that is,
\[
    |{\det}_2(I+B_\kappa)| \int_{\mathcal{W}}
     \phi(\iota+F_\kappa) 
     e^{\mathfrak{q}_{\eta(\kappa)}} d\mu
    =e^{\|\kappa\|_2^2/2} \int_{\mathcal{W}} \phi d\mu.
\]
Substituting $\phi=f(\iota+F_{\widehat{\kappa}})$ and
\eqref{t.inv.transf.1} into this identity,
we obtain \eqref{eq:inv.transf}.
\end{proof}

Looking at the identity \eqref{eq:inv.transf} from the point
of view of $\widehat{\kappa}$, and switching $\kappa$ and
$\widehat{\kappa}$, we obtain the following.   

\begin{corollary}\label{c.inv.transf}
Let $\kappa\in \mathcal{L}_2$ and assume that 
${\det}_2(I+B_\kappa)\ne0$.
Then $\Lambda(B_{\eta(\widehat{\kappa})})<1$ and 
the distribution
$\mu\circ(\iota+F_\kappa)^{-1}$ of 
$\iota+F_\kappa$ under $\mu$ has the Radon-Nikodym
derivative with respect to $\mu$ of the form  
$|{\det}_2(I+B_{\widehat{\kappa}})|
     e^{-\|\widehat{\kappa}\|_2^2/2} 
     e^{\mathfrak{q}_{\eta(\widehat{\kappa})}}$, i.e., 
\[
    \frac{d(\mu\circ(\iota+F_\kappa)^{-1})}{d\mu}
    =|{\det}_2(I+B_{\widehat{\kappa}})|
     e^{-\|\widehat{\kappa}\|_2^2/2} 
     e^{\mathfrak{q}_{\eta(\widehat{\kappa})}}.
\]
\end{corollary}

\begin{proof}
Due to the identity 
$(I+B_\kappa)(I+B_{\widehat{\kappa}})=I$, it holds that
\[
    \|h\|_{\mathcal{H}}
    \le \|(I+B_\kappa)\|_{\text{\rm op}}
       \|(I+B_{\widehat{\kappa}})h\|_{\mathcal{H}}
    \quad\text{for }h\in \mathcal{H}.
\]
Hence 
\[
    \inf_{\|h\|_{\mathcal{H}}=1}
      \|(I+B_{\widehat{\kappa}})h\|_{\mathcal{H}}^2
     \ge \|(I+B_\kappa)\|_{\text{\rm op}}^{-2}>0.
\]
As was mentioned in Remark~\ref{r.transf}(i), this implies
that $\Lambda(B_{\eta(\widehat{\kappa})})<1$.

If we set $\kappa^\prime=\widehat{\kappa}$, then 
$\widehat{\kappa^\prime}=\kappa$ in $\mathcal{L}_2$.
Applying Theorem~\ref{t.inv.transf} to 
$\widehat{\kappa}$, we see that
\[
    |{\det}_2(I+B_{\widehat{\kappa}})|
     e^{-\|\widehat{\kappa}\|_2^2/2}
    \int_{\mathcal{W}} f
      e^{\mathfrak{q}_{\eta(\widehat{\kappa})}}
      d\mu
    =\int_{\mathcal{W}} 
       f(\iota+F_\kappa)d\mu
\]
for every $f\in C_b(\mathcal{W})$.
Thus the desired expression of the Radon-Nikodym derivative
is obtained.
\end{proof}

\section{Surjectivity and injectivity}
\label{sec:surj}
In this section, we present a method to find
$\kappa\in\mathcal{L}_2$ satisfying that
$\eta(\kappa)=\eta$ in $\mathcal{L}_2$ for given 
$\eta\in\mathcal{S}_2$ with
$\Lambda(B_\eta)<1$, and show the relating ``bijectivity''
of the mapping $\kappa\mapsto \eta(\kappa)$.
To state the method, let 
$\mathcal{S}_+(\mathcal{H})$ be the totality of all 
self-adjoint continuous linear operators 
$C:\mathcal{H}\to \mathcal{H}$ with $C\ge0$, that is,
$\langle Ch,h\rangle_{\mathcal{H}}\ge0$ for every 
$h\in \mathcal{H}$.
Given $\eta\in\mathcal{S}_2$ with $\Lambda(B_\eta)<1$, 
it holds that 
$I-B_\eta\in \mathcal{S}_+(\mathcal{H})$.
Due to the square root lemma (\cite{RN}) asserting
that every $A\in \mathcal{S}_+(\mathcal{H})$ has a unique  
$C\in \mathcal{S}_+(\mathcal{H})$ with $C^2=A$, there is a
unique $C_\eta\in \mathcal{S}_+(\mathcal{H})$ with
$C_\eta^2=I-B_\eta$.
Since $C_\eta\ge0$, 
$\|(I+C_\eta)h\|_{\mathcal{H}}^2\ge\|h\|_{\mathcal{H}}^2$
for $h\in \mathcal{H}$.
Hence $I+C_\eta$ has a continuous inverse
(cf.~the argument in Remark~\ref{r.inv.transf}(ii)).
Using this inverse, we see that
$C_\eta-I=-(I+C_\eta)^{-1}B_\eta
 \in \mathcal{S}(\mathcal{H}^{\otimes 2})$.
By Lemma~\ref{l.b.kappa}, there is a 
$\kappa_S(\eta)\in \mathcal{S}_2$ with
$B_{\kappa_S(\eta)}=C_\eta-I$.
The following is the method to find $\kappa$ with 
$\eta(\kappa)=\eta$ in $\mathcal{L}_2$.

\begin{theorem}\label{t.surj.1st}
Let $\eta\in\mathcal{S}_2$ and assume that
$\Lambda(B_\eta)<1$. 
Then the following assertions hold.
\\
{\rm (i)} If $\kappa=\kappa_S(\eta)$, then 
$\eta(\kappa)=\eta$ in $\mathcal{L}_2$.
\\
{\rm (ii)} It holds that
\begin{equation}\label{t.surj.1st.1}
    \int_{\mathcal{W}} fe^{\mathfrak{q}_\eta} d\mu
    =\{{\det}_2(I-B_\eta)\}^{-1/2} 
     \int_{\mathcal{W}} f(\iota+F_{\widehat{\kappa_S(\eta)}})d\mu
\end{equation}
for every $f\in C_b(\mathcal{W})$.
\end{theorem}

\begin{remark}\label{r.surj.1st}
(i) By the first assertion,
$\Lambda(B_{\eta(\kappa_S(\eta))})<1$.
As was seen before Theorem~\ref{t.inv.transf}, 
$\widehat{\kappa_S(\eta)}$ is defined well.
\\
(ii) For $\eta\in\mathcal{S}_2$, the Laplace transformation
of the distribution of $\mathfrak{q}_\eta$ under the signed
measure $fd\mu$ is the function 
$\mathbb{R}\ni \lambda\mapsto
 \int_{\mathcal{W}} f e^{\lambda \mathfrak{q}_\eta}d\mu$.
Since $\lambda\mathfrak{q}_\eta=\mathfrak{q}_{\lambda\eta}$, 
\eqref{t.surj.1st.1} gives an evaluation of the Laplace
transformation.
\\
(iii)
An identity similar to \eqref{t.surj.1st.1} was shown in
\cite[Theorem\,5.7.6]{mt-cambridge} under the much stronger
assumption that $\|B_\eta\|_{\text{\rm op}}<3/8$ than the
one that $\Lambda(B_\eta)<1$
(cf.~Remark~\ref{r.transf}(i)).
Therein $F_{\widehat{\kappa_S(\eta)}}$ was constructed as a
series in $L^2(\mu;\mathcal{H})$ like the one in
Lemma~\ref{l.f.kappa} by using the eigenfunction expansion
of $B_\eta$. 
\\
(iv)
As for the Fourier transformation, i.e., the function
$\mathbb{R}\ni \lambda\mapsto
 \int_{\mathcal{W}} 
  fe^{\sqrt{-1}\lambda\mathfrak{q}_\eta} d\mu$,
a formula similar to
\eqref{t.surj.1st.1} was shown by Malliavin and the author
\cite{mall-t} for ``analytic'' $f$, by complexifying
$F_{\widehat{\kappa_S(\eta)}}$ with the help of the
eigenfunction expansion of $\sqrt{-1}B_\eta$. 
\end{remark}

\begin{proof}
Let $C_\eta$ be as above and 
$\kappa=\kappa_S(\eta)$.
\\
(i)
Since 
$B_\kappa\in \mathcal{S}(\mathcal{H}^{\otimes 2})$,
by \eqref{r.transf.1} and the definition of
$\kappa_S(\eta)$, we have that
\[
    -B_{\eta(\kappa)}=(I+B_\kappa)^2-I=-B_\eta.
\]
Due to Lemmma~\ref{l.b.kappa}, $\eta(\kappa)=\eta$
in $\mathcal{L}_2$.
\\
(ii)
It follows from \eqref{r.inv.transf.1} that
\[
    {\det}_2\bigl((I+B)^2\bigr)
    =\{{\det}_2(I+B)\}^2 
      \exp(-\|B\|_{\mathcal{H}^{\otimes 2}}^2)
    \quad\text{for }
    B\in \mathcal{S}(\mathcal{H}^{\otimes 2}).
\]
Since $\|B_\kappa\|_{\mathcal{H}^{\otimes 2}}=\|\kappa\|_2$,
this identity implies that
\[
    |{\det}_2(I+B_\kappa)|e^{-\|\kappa\|_2^2/2}
    =\{{\det}_2((I+B_\kappa)^2)\}^{1/2}
    =\{{\det}_2(I-B_\eta)\}^{1/2}.
\]
Substituting this and the equality 
$\mathfrak{q}_{\eta(\kappa)}=\mathfrak{q}_\eta$ into
the second identity in Theorem~\ref{t.inv.transf},
we obtain \eqref{t.surj.1st.1}.
\end{proof}

As was seen in Remark~\ref{r.transf}(ii), the mapping
$\mathcal{L}_2\ni\kappa
 \mapsto\eta(\kappa)\in\mathcal{S}_2$ 
is not injective. 
Restricting the domain of the mapping, we have the following 
bijectivity.

\begin{proposition}\label{p.inj.1st}
Put
$\mathcal{S}_{2,+}=\bigl\{\kappa\in\mathcal{S}_2;
      I+B_\kappa\ge0\bigr\}$,
$\widehat{\mathcal{S}_{2,+}}
 =\bigl\{\kappa\in \mathcal{S}_{2,+};
  \Lambda(B_{\eta(\kappa)})<1\}$,
and
$\widehat{\mathcal{S}_2}
 =\{\eta\in\mathcal{S}_2;
     \Lambda(B_\eta)<1\}$.
\\
{\rm(i)}
If $\kappa_1,\kappa_2\in \mathcal{S}_{2,+}$ and 
$\eta(\kappa_1)=\eta(\kappa_2)$ in $\mathcal{L}_2$, then  
$\kappa_1=\kappa_2$ in $\mathcal{L}_2$.
\\
{\rm(ii)}
The mapping 
$\widehat{\mathcal{S}_{2,+}}\ni\kappa\mapsto
 \eta(\kappa)\in\widehat{\mathcal{S}_2}$ 
is bijective in the sense that
$\kappa_1=\kappa_2$ in $\mathcal{L}_2$ if 
$\kappa_1,\kappa_2\in\widehat{\mathcal{S}_{2,+}}$ and
$\eta(\kappa_1)=\eta(\kappa_2)$ in $\mathcal{L}_2$,
and each $\eta\in\widehat{\mathcal{S}_2}$ admits a
$\kappa\in \widehat{\mathcal{S}_{2,+}}$ with
$\eta=\eta(\kappa)$ in $\mathcal{L}_2$.
\end{proposition}

\begin{proof}
The assertion~(ii) is an immediate consequence of the
assertion~(i) and Theorem~\ref{t.surj.1st}.
To see the assertion~(i), let 
$\kappa_1,\kappa_2\in \mathcal{S}_{2,+}$ and assume that  
$\eta(\kappa_1)=\eta(\kappa_2)$ in $\mathcal{L}_2$.
By \eqref{r.transf.1}, 
$(I+B_{\kappa_1})^2=(I+B_{\kappa_2})^2$.
Since 
$I+B_{\kappa_1},I+B_{\kappa_2}
 \in \mathcal{S}_+(\mathcal{H})$, due to the square root
lemma, $I+B_{\kappa_1}=I+B_{\kappa_2}$. 
By Lemma~\ref{l.b.kappa}, $\kappa_1=\kappa_2$ in
$\mathcal{L}_2$. 
\end{proof}

In the remaining of this section, we give an application of
Theorem~\ref{t.surj.1st}. 
For $\kappa\in\mathcal{L}_2$ and $x\in \mathbb{R}^d$, define 
$\mathfrak{h}(\kappa;x),\mathfrak{h}(\kappa):
 \mathcal{W}\to \mathbb{R}$ by 
\[
    \mathfrak{h}(\kappa;x)
      =\frac12\int_0^T\biggl\langle x,
        \int_0^T\kappa(t,s)d\theta(s)
       \biggr\rangle^2 dt
       \quad\text{and}\quad
     \mathfrak{h}(\kappa)
      =\frac12\int_0^T\biggl|\int_0^T
        \kappa(t,s)d\theta(s)\biggr|^2 dt.
\]
If $d=1$ and 
$\kappa(t,s)=\boldsymbol{1}_{[0,t)}(s)$ for 
$(t,s)\in[0,T]^2$,
then
$\mathfrak{h}(\kappa;1)=\mathfrak{h}(\kappa)
 =\int_0^T\theta(t)^2 dt/2$, 
which relates to the harmonic oscillator 
$-(1/2)\{(d/dx)^2-x^2\}$,
one of the fundamental Schr\"odinger operators.
In the stochastic approach to the KdV equation,
Ikeda and the author \cite{ikeda-t,st-spa} used
$\mathfrak{h}(\kappa;x)$ with $\kappa$ of the form 
$\kappa(t,s)
    =\boldsymbol{1}_{[0,t)}(s)
     \text{\rm diag}[
        e^{(t-s)p_1},\dots, e^{(t-s)p_d}]$
for $(t,s)\in[0,T]^2$, 
where $p_1,\dots,p_d\in \mathbb{R}$ and $\text{\rm
  diag}[a_1,\dots,a_d]$ is the $d$-dimensional diagonal
matrix whose $i$th entry is $a_i$.

\begin{proposition}\label{p.h.osc}
Let $\kappa\in\mathcal{L}_2$ and $x\in \mathbb{R}^d$.
Define 
$c(\kappa;x),c(\kappa),
 c^\prime(\kappa;x),c^\prime(\kappa)
 \in\mathcal{S}_2$ 
by
\begin{align*}
    c(\kappa;x)(t,s)
    & =\int_0^T \bigl[(\kappa(u,s)^\dagger x)
          \otimes (\kappa(u,t)^\dagger x)\bigr] du,
    \\
    c(\kappa)(t,s)
    & =\int_0^T \kappa(u,t)^\dagger \kappa(u,s)du
    \quad \text{for } (t,s)\in[0,T]^2
\end{align*}
and $c^\prime(\kappa;x)=\kappa_S(-c(\kappa;x))$,
$c^\prime(\kappa)=\kappa_S(-c(\kappa))$.
Then 
$\Lambda(B_{-c(\kappa;x)})\le0$,
$\Lambda(B_{-c(\kappa)})\le0$, 
$B_{c(\kappa;x)}$ is of trace class, and it holds that 
\begin{align}
    \int_{\mathcal{W}} f e^{-\mathfrak{h}(\kappa;x)}
         d\mu
    & =\{\det(I+B_{c(\kappa;x)})\}^{-1/2}
     \int_{\mathcal{W}} f(\iota
         +F_{\widehat{c^\prime(\kappa;x)}})d\mu
 \label{p.h.osc.1}
    \\
    \int_{\mathcal{W}} f e^{-\mathfrak{h}(\kappa)} d\mu
    & =\{\det(I+B_\kappa^*B_\kappa)\}^{-1/2}
     \int_{\mathcal{W}} f(\iota
         +F_{\widehat{c^\prime(\kappa)}})d\mu
 \label{p.h.osc.2}
\end{align}
for every $f\in C_b(\mathcal{W})$.
\end{proposition}

\begin{proof}
Observe that 
\begin{align*}
    \langle D\mathfrak{h}(\kappa;x),g
      \rangle_{\mathcal{H}}
    &  = \int_0^T \biggl\langle x,\int_0^T \kappa(t,s)
           d\theta(s) \biggr\rangle
         \biggl\langle x,\int_0^T \kappa(t,u)
         g^\prime(u) du \biggr\rangle dt,
   \\
    \langle D^2\mathfrak{h}(\kappa;x),h\otimes g
      \rangle_{\mathcal{H}^{\otimes 2}}
    & = \int_0^T \biggl\langle x,\int_0^T \kappa(t,s)
         h^\prime(s)ds \biggr\rangle
        \biggl\langle x,\int_0^T \kappa(t,u)
         g^\prime(u)du \biggr\rangle dt
   \\
    & =\langle B_{c(\kappa;x)},h\otimes g
        \rangle_{\mathcal{H}^{\otimes 2}}
    \quad\text{for }h,g\in \mathcal{H}.
\end{align*}
Hence 
\[
    \int_{\mathcal{W}}D \mathfrak{h}(\kappa;x) d\mu=0
    \quad\text{and}\quad
    D^2 \mathfrak{h}(\kappa;x)=B_{c(\kappa;x)}.
\]
Rewriting 
$\bigl\langle x,\int_0^T\kappa(t,s)d\theta(s)
 \bigr\rangle$
as 
$\int_0^T\langle\kappa(t,s)^\dagger x,d\theta(s)\rangle$, 
we see that
\[
    \int_{\mathcal{W}} \mathfrak{h}(\kappa;x)d\mu
    =\frac12\int_0^T\int_0^T|\kappa(t,s)^\dagger x|^2 
     dsdt.
\]
By \cite[Proposition~5.7.4]{mt-cambridge} and
Lemma~\ref{l.q.eta}, it holds that
\begin{equation}\label{p.h.osc.21}
    -\mathfrak{h}(\kappa;x)
    =\mathfrak{q}_{-c(\kappa;x)}
     -\frac12 \int_0^T\int_0^T|\kappa(t,s)^\dagger x|^2 
       dsdt.
\end{equation}

Since
\begin{equation}\label{p.h.osc.22}
    \langle B_{-c(\kappa;x)}h,h\rangle_{\mathcal{H}}
    =-\int_0^T\biggl(\int_0^T \langle 
       \kappa(t,s)^\dagger x,h^\prime(s)\rangle ds
       \biggr)^2 dt
    \quad\text{for } h\in \mathcal{H},
\end{equation}
we see that
\[
    \Lambda(B_{-c(\kappa;x)})\le 0.
\]
By Theorem~\ref{t.surj.1st} and \eqref{p.h.osc.21}, we
have that
\begin{align}
    & \int_{\mathcal{W}} f e^{-\mathfrak{h}(\kappa;x)} d\mu
 \nonumber
    \\
    & \quad
     =\exp\biggl(-\frac12\int_0^T\int_0^T
           |\kappa(t,s)^\dagger x|^2 dsdt\biggr)
     \{{\det}_2(I+B_{c(\kappa;x)})\}^{-1/2}
     \int_{\mathcal{W}}
       f(\iota+F_{\widehat{c^\prime(\kappa;x)}}) d\mu
\label{p.h.osc.23}
\end{align}
for every $f\in C_b(\mathcal{W})$.

For $A\in \mathcal{H}^{\otimes 2}$, denote by 
$|A|\in \mathcal{S}_+(\mathcal{H})$ the square root of
$A^*A$.
By definition, $A$ is of trace class if 
$\sum_{n=1}^\infty \langle |A|h_n,h_n
  \rangle_{\mathcal{H}}<\infty$
for some orthonormal basis $\{h_n\}_{n=1}^\infty$
of $\mathcal{H}$.
Since 
$B_{c(\kappa;x)}\in \mathcal{S}(\mathcal{H}^{\otimes 2})$,
$B_{c(\kappa;x)}^*B_{c(\kappa;x)}=B_{c(\kappa;x)}^2$.
By \eqref{p.h.osc.22}, 
$B_{c(\kappa;x)}\in \mathcal{S}_+(\mathcal{H})$.
Due to the uniqueness of square root, 
$B_{c(\kappa;x)}=|B_{c(\kappa;x)}|$.
By \eqref{p.h.osc.22} again, we have that
\[
    \sum_{n=1}^\infty \langle B_{c(\kappa;x)}h_n,h_n
       \rangle_{\mathcal{H}}
    =\int_0^T\int_0^T|\kappa(t,s)^\dagger x|^2 dsdt
\]
for any orthonormal basis $\{h_n\}_{n=1}^\infty$ of
$\mathcal{H}$.
Hence $B_{c(\kappa;x)}$ is of trace class and 
\[
    \text{\rm tr} B_{c(\kappa;x)}
    =\int_0^T\int_0^T|\kappa(t,s)^\dagger x|^2 dsdt.
\]
Thus we have that 
\[
    {\det}_2(I+B_{c(\kappa;x)})
    =\det(I+B_{c(\kappa;x)})
     \exp\biggl(-\int_0^T\int_0^T|\kappa(t,s)^\dagger x|^2 
       dsdt\biggr).
\]
Plugging this into \eqref{p.h.osc.23}, we obtain
\eqref{p.h.osc.1}. 

Let $e_1,\dots,e_d$ be an orthonormal basis of
$\mathbb{R}^d$.
Notice that    
$\mathfrak{h}(\kappa)
    =\sum_{i=1}^d \mathfrak{h}(\kappa;e_i)$.
Since 
\[
    \sum_{i=1}^d c(\kappa;e_i)=c(\kappa)
    ~~\text{and}~~
    \sum_{i=1}^d \int_0^T\int_0^T
      |\kappa(t,s)^\dagger e_i|^2 dsdt=\|\kappa\|_2^2,
\]
by \eqref{p.h.osc.21}, we have that
\[
    -\mathfrak{h}(\kappa)
    =\mathfrak{q}_{-c(\kappa)}-\frac12 \|\kappa\|_2^2.
\]
Further, by \eqref{p.h.osc.22},
\[
    \Lambda(B_{-c(\kappa)})\le 0. 
\]
Due to Theorem~\ref{t.surj.1st}, it holds that
\[
    \int_{\mathcal{W}} f e^{-\mathfrak{h}(\kappa)} d\mu
    =e^{-\|\kappa\|_2^2/2}
     \{{\det}_2(I+B_{c(\kappa)})\}^{-1/2}
     \int_{\mathcal{W}} f(\iota+F_{\widehat{c^\prime(\kappa)}})
     d\mu
\]
for every $f\in C_b(\mathcal{W})$.
Since $B_{c(\kappa)}=B_\kappa^*B_\kappa$, $B_{c(\kappa)}$ is
of trace class and  
$\text{\rm tr} B_{c(\kappa)}
 =\|B_\kappa\|_{\mathcal{H}^{\otimes 2}}^2=\|\kappa\|_2^2$.
Plugging these into the above identity, we obtain 
\eqref{p.h.osc.2}.
\end{proof}

\section{Linear transformations}
\label{sec:lin.transf}

In this section, we apply Theorem~\ref{t.transf} to linear
transformations studied by Cameron and Martin \cite{CM}.
Precisely speaking, define
$\mathbb{F}_\phi:\mathcal{W}\to \mathcal{H}$ 
for $\phi\in\mathcal{L}_2$ by
\[
    \langle \mathbb{F}_\phi, h\rangle_{\mathcal{H}}
    =\int_0^T \biggl\langle 
       \int_0^T \phi(t,s)\theta(s)ds, h^\prime(t)
       \biggr\rangle dt
    \quad\text{for } h\in \mathcal{H}.
\]
We investigate the linear transformation
\[
    \iota+\mathbb{F}_\phi:\mathcal{W}\in w \mapsto
     w+\int_0^\bullet\biggl(\int_0^T \phi(t,s)w(s)ds
       \biggr)dt \in \mathcal{W}.
\]
Set $\kappa_\phi\in\mathcal{L}_2$ so that 
\[
    \kappa_\phi(t,s)=\int_s^T \phi(t,u)du
    \quad\text{for } (t,s)\in[0,T]^2.
\]
Applying It\^o's formula, we see that
\begin{equation}\label{eq:kappa.phi}
    \int_0^T \kappa_\phi(t,s)d\theta(s)
    =\int_0^T \phi(t,s)\theta(s)ds
    \quad\text{for } t\in[0,T].
\end{equation}
Hence $\mathbb{F}_\phi=F_{\kappa_\phi}$.
The aim of this section is to show the following theorem.
The comparison between our result and Cameron and Martin's
will be given in Remark~\ref{r.transf.cm} below.

\begin{theorem}\label{t.transf.cm}
Let $\phi\in\mathcal{L}_2$ and assume that
$\Lambda(B_{\eta(\kappa_\phi)})<1$. 
Put
\[
    \Psi_\phi 
    = -\int_0^T\biggl\langle \int_0^T 
         \phi(t,s)^\dagger d\theta(t),
       \theta(s)\biggr\rangle ds
      -\frac12 \int_0^T\biggl|\int_0^T
         \phi(t,s)\theta(s)ds \biggr|^2dt.
\]
Then $B_{\kappa_\phi}$ is of trace class and it holds that
\begin{equation}\label{t.transf.cm.1}
    |\det(I+B_{\kappa_\phi})|\int_{\mathcal{W}} 
      f(\iota+\mathbb{F}_\phi) e^{\Psi_\phi}
      d\mu
    =\int_{\mathcal{W}} f d\mu
\end{equation}
for every $f\in C_b(\mathcal{W})$.
\end{theorem}

For the proof, we prepare lemmas.

\begin{lemma}\label{l.d*f.kappa}
It holds that 
\[
    \mathfrak{q}_{\eta(\kappa)}
     =-D^*F_\kappa-\frac12 \|F_\kappa\|_{\mathcal{H}}^2
    +\frac12\|\kappa\|_2^2
\]
for every $\kappa\in\mathcal{L}_2$.
\end{lemma}

\begin{proof}
Let $n\in \mathbb{N}$, $h_1,\dots,h_n\in \mathcal{H}$ and 
$\varphi\in C_b^\infty(\mathbb{R}^n)$($\equiv$ the space of 
real $C^\infty$-functions on $\mathbb{R}^n$ whose
derivatives of all orders are bounded). 
Put 
$D^*\boldsymbol{h}=(D^*h_1,\dots,D^*h_n)
 \in \mathbb{D}^\infty(\mathbb{R}^n)$.
By Lemma~\ref{l.f.kappa}, the chain rule for $D$ 
(\cite[Corollay~5.3.2]{mt-cambridge}), and the symmetry of
the Hessian matrix  
$\bigl(\partial^2\varphi/\partial x^i \partial x^j
 \bigr)_{1\le i,j\le n}$, we see that
\begin{align*}
    \int_{\mathcal{W}} (D^*F_\kappa)
          \varphi(D^*\boldsymbol{h}) d\mu
    & =\int_{\mathcal{W}} \langle B_\kappa,
          D^2\bigl(\varphi(D^*\boldsymbol{h})\bigr) 
          \rangle_{\mathcal{H}^{\otimes 2}} d\mu
    \\
    & =\sum_{i,j=1}^n \int_{\mathcal{W}} \langle B_\kappa,
         h_i\otimes h_j\rangle_{\mathcal{H}^{\otimes 2}} 
         \frac{\partial^2\varphi}{
                \partial x^i \partial x^j}(
           D^*\boldsymbol{h}) d\mu
    \\
    & =\frac12 \sum_{i,j=1}^n \int_{\mathcal{W}} 
        \langle (B_\kappa+B_\kappa^*),
         h_i\otimes h_j\rangle_{\mathcal{H}^{\otimes 2}} 
         \frac{\partial^2\varphi}{
                \partial x^i \partial x^j}(
           D^*\boldsymbol{h}) d\mu
    \\
    & =\frac12
         \int_{\mathcal{W}} 
           \bigl((D^*)^2(B_\kappa+B_\kappa^*)\bigr)
           \varphi(D^*\boldsymbol{h}) d\mu.
\end{align*}
Thus we obtain that 
\begin{equation}\label{l.d*f.kappa.21}
    D^*F_\kappa=\frac12(D^*)^2(B_\kappa+B_\kappa^*).
\end{equation}

Since $D F_\kappa=B_\kappa$ as was seen in the proof of
Lemma~\ref{l.f.kappa}, we have that
\begin{align*}
     \biggl\langle D\biggl(\frac12 
       \|F_\kappa\|_{\mathcal{H}}^2\biggr),g
      \biggr\rangle_{\mathcal{H}}
    &  =\langle F_\kappa,B_\kappa g\rangle_{\mathcal{H}},
    \\
     \biggl\langle D^2\biggl(\frac12
         \|F_\kappa\|_{\mathcal{H}}^2\biggr),
        h\otimes g \biggr\rangle_{\mathcal{H}^{\otimes 2}}
    &  =\langle B_\kappa h,B_\kappa g\rangle_{\mathcal{H}}
      =\langle B_\kappa^*B_\kappa, h\otimes g
       \rangle_{\mathcal{H}^{\otimes 2}}
     \quad \text{for }h,g\in \mathcal{H}.
\end{align*}
Hence
\[
    \int_{\mathcal{W}} D\biggl(
       \frac12\|F_\kappa\|_{\mathcal{H}}^2\biggr) d\mu=0
    \quad\text{and}\quad
    D^2\biggl(\frac12\|F_\kappa\|_{\mathcal{H}}^2\biggr)
    =B_\kappa^*B_\kappa. 
\]
Further, by the It\^o isometry, we see that
\[
    \int_{\mathcal{W}}\biggl(\frac12
        \|F_\kappa\|_{\mathcal{H}}^2\biggr) d\mu
    =\frac12\|\kappa\|_2^2.
\]
Then, by \cite[Proposition~5.7.4]{mt-cambridge}, we obtain 
that 
\begin{equation}\label{l.d*f.kappa.22}
    \frac12 \|F_\kappa\|_{\mathcal{H}}^2
    =\frac12(D^*)^2(B_\kappa^*B_\kappa)
     +\frac12 \|\kappa\|_2^2.
\end{equation}
Applying Lemma~\ref{l.q.eta} to 
$\mathfrak{q}_{\eta(\kappa)}$
with the help of \eqref{r.transf.1}, we know that
\[
    \mathfrak{q}_{\eta(\kappa)}
    =-\frac12(D^*)^2\{B_\kappa+B_\kappa^*
     +B_\kappa^*B_\kappa\}.
\]
Plugging \eqref{l.d*f.kappa.21} and \eqref{l.d*f.kappa.22}
into this equality, we obtain the desired identity.
\end{proof}

\begin{lemma}\label{l.transf.cm}
Let $\phi\in\mathcal{L}_2$ and assume that
$\Lambda(B_{\eta(\kappa_\phi)})<1$.
Put 
\[
    \widetilde{\Psi}_\phi 
    = \Psi_\phi
    +\int_0^T\biggl(\int_0^s 
       [\text{\rm tr}\phi(t,s)]dt\biggr)ds.
\]
Then it holds that 
\begin{equation}\label{l.transf.cm.1}
    |{\det}_2(I+B_{\kappa_\phi})|\int_{\mathcal{W}} 
      f(\iota+\mathbb{F}_\phi) 
      e^{\widetilde{\Psi}_\phi} d\mu
    =\int_{\mathcal{W}} f d\mu
\end{equation}
for every $f\in C_b(\mathcal{W})$.
\end{lemma}

\begin{proof}
By Theorem~\ref{t.transf} and Lemma~\ref{l.d*f.kappa},
it suffices to show the identity
\begin{equation}\label{l.transf.cm.21}
    D^*F_{\kappa_\phi}=\int_0^T\biggl\langle\int_0^T
       \phi(t,s)^\dagger d\theta(t),\theta(s)
       \biggr\rangle ds
    -\int_0^T\biggl(\int_0^s 
        [\text{\rm tr}\phi(t,s)]dt\biggr)ds.
\end{equation}
Let $\phi_m,\phi\in\mathcal{L}_2$, $m\in \mathbb{N}$, and
suppose that 
$\|\phi_m-\phi\|_2\to0$ as $m\to\infty$.
Then 
\[
    \|B_{\kappa_{\phi_m}}-B_{\kappa_\phi}
       \|_{\mathcal{H}^{\otimes 2}}
    =\|\kappa_{\phi_m}-\kappa_\phi\|_2\to0.
\]
By Lemma~\ref{l.f.kappa} and the continuity of $D^*$, 
$F_{\kappa_{\phi_m}}\to F_{\kappa_\phi}$ in
$L^p(\mu;\mathcal{H})$ for any $p\in(1,\infty)$.
Further it is easily seen that the right hand side of
\eqref{l.transf.cm.21} for $\phi_m$ converges to that for
$\phi$ in probability. 
Thus, to show \eqref{l.transf.cm.21}, we may and will assume
that $\phi\in\mathcal{L}_2$ is of the form 
\[
    \phi(t,s)=\sum_{n=0}^{N-1} 
       \boldsymbol{1}_{[T(N;n),T(N;n+1))}(t)
       \phi(T(N;n),s)
    \quad\text{for } (t,s)\in[0,T]^2
\]
with some $N\in \mathbb{N}$, where $T(N;n)=nT/N$.
Let $e_1,\dots,e_d$ be an orthonormal basis of
$\mathbb{R}^d$.
Define 
$h_{n;i},k_{s;i}\in \mathcal{H}$ by 
$h_{n;i}^\prime=\boldsymbol{1}_{[T(N;n),T(N;n+1))}e_i$ and
$k_{s;i}^\prime=\boldsymbol{1}_{[0,s)}e_i$
for $1\le i\le d$, $0\le n\le N-1$, and $s\in[0,T]$.
By \eqref{eq:kappa.phi}, it holds that
\[
    F_{\kappa_\phi}=\sum_{n=0}^{N-1} \sum_{i,j=1}^d
    \biggl(\int_0^T \phi_j^i(T(N;n),s)
    \theta^j(s)ds\biggr) h_{n;i}.
\]
Since 
$D^*h_{n;i}=\theta^i(T(N;n+1))-\theta^i(T(N;n))$ and 
$D\theta^i(s)=k_{s;i}$, due to the product rule for $D^*$ 
(\cite[Theorem~5.2.8]{mt-cambridge}), we obtain that
\begin{align*}
    D^*F_{\kappa_\phi}
    =& \sum_{n=0}^{N-1} \sum_{i,j=1}^d
    \biggl(\int_0^T \phi_j^i(T(N;n),s)
    \theta^j(s)ds\biggr)
    \{\theta^i(T(N;n+1))-\theta^i(T(N;n))\}
    \\
    & -\sum_{n=0}^{N-1} \sum_{i,j=1}^d
       \biggl\langle \int_0^T 
       \phi_j^i(T(N;n),s) k_{s;j} ds,
       h_{n;i}\biggr\rangle_{\mathcal{H}}.
\end{align*}
The first term of the right hand side of the equality
coincides with 
\begin{align*}
    & \int_0^T \biggl[\sum_{j=1}^d \biggl(\sum_{i=1}^d
       \biggl\{\sum_{n=0}^{N-1}
         \phi_j^i(T(N;n),s)
         \{\theta^i(T(N;n+1)-\theta^i(T(N;n))\}
       \biggr\}\biggr)\theta^j(s)\biggr]ds
    \\
    & =\int_0^T \biggl\langle \int_0^T \phi(t,s)^\dagger
        d\theta(t),\theta(s)\biggr\rangle ds.
\end{align*}
Since 
\[
    \langle k_{s;j},h_{n;i}\rangle_{\mathcal{H}}
    =\bigl\{(T(N;n+1)\wedge s)
         -(T(N;n)\wedge s)\bigr\}\delta_{ji},
\]
the second term of the right hand side coincides with
\begin{align*}
     & \int_0^T \biggl[\sum_{n=0}^{N-1} \biggl(
         \sum_{i=1}^d
         \phi_i^i(T(N;n),s)\biggr)
         \bigl\{(T(N;n+1)\wedge s)
           -(T(N;n)\wedge s)\bigr\}\biggr] ds
     \\
     & =\int_0^T \biggl(\int_0^s [\text{\rm tr}\phi(t,s)]dt
      \biggr)ds.
\end{align*}
Thus \eqref{l.transf.cm.21} holds.
\end{proof}

\begin{lemma}\label{l.kappa.phi.trace}
$B_{\kappa_\phi}$ is of trace class.
\end{lemma}

\begin{proof}
Define $\psi\in\mathcal{L}_2$ by 
$\psi(t,s)=\boldsymbol{1}_{[0,t)}(s)I_d$
for $(t,s)\in[0,T]^2$.
It holds that 
\[
    \kappa_\phi(t,s)
    =\int_0^T \phi(t,u)\psi(u,s)du
    \quad\text{for }(t,s)\in[0,T]^2.
\]
Hence
\begin{equation}\label{l.kappa.phi.trace.21}
    B_{\kappa_\phi}=B_\phi B_\psi,
\end{equation}
which implies that $B_{\kappa_\phi}$ is of trace class.
\end{proof}

\begin{lemma}\label{l.trace}
If $\kappa\in\mathcal{L}_2$ is continuous on $[0,T]^2$ and 
$B_\kappa$ is of trace class, then 
\[
    \text{\rm tr} B_\kappa
    =\int_0^T [\text{\rm tr}\kappa(s,s)]ds.
\]
\end{lemma}
The continuity of $\kappa$ is used to determine
$\kappa(s,s)$ for each $s\in[0,T]$.

\begin{proof}
The assertion can be shown as an extension of the proof of 
\cite[Theorem~IV.8.1]{GGK} to general dimensions, which
theorem deals with the one-dimensional case.
The extension is based on the fact that the
Cameron-Martin subspace is the $d$-times product space of   
one-dimensional ones, and it is routine.
We omit the details.
\end{proof}

\begin{proof}[Proof of Theorem~\ref{t.transf.cm}]
By Lemma~\ref{l.kappa.phi.trace}, it suffices to show
\eqref{t.transf.cm.1}.
Suppose that $\phi_m\in\mathcal{L}_2$, $m\in \mathbb{N}$,
and $\|\phi_m-\phi\|_2\to0$ as $m\to\infty$.
Denote by $\|\cdot\|_1$ the trace class norm.
By \cite[Lemma~IV.7.2]{GGK} and
\eqref{l.kappa.phi.trace.21}, we have that  
\[
    \|B_{\kappa_{\phi_m}}-B_{\kappa_\phi}\|_1
    \le \|B_\psi\|_{\mathcal{H}^{\otimes 2}}
        \|B_{\phi_m}-B_\phi\|_{\mathcal{H}^{\otimes 2}}
    =\|B_\psi\|_{\mathcal{H}^{\otimes 2}}
     \|\phi_m-\phi\|_2
    \to0 \quad\text{as }m\to\infty.
\]
Thanks to the continuity of the mapping $B\mapsto\det(I+B)$
with respect to the trace class norm
(\cite[Corollary~II.4.2]{GGK}), 
$\det(I+B_{\kappa_{\phi_m}})$ converges to
$\det(I+B_{\kappa_\phi})$.
It is easily seen that $\Psi_{\phi_m}$ converges to
$\Psi_\phi$ in probability.
Further,
$|\Lambda(B_{\kappa_{\phi_m}})
  -\Lambda(B_{\kappa_\phi})|
 \le \|B_{\kappa_{\phi_m}}
   -B_{\kappa_\phi}\|_{\mathcal{H}^{\otimes 2}}\to0$ 
as $m\to\infty$, and hence 
$\sup_{m\ge m_0}\Lambda(B_{\kappa_{\phi_m}})<1$ for some
$m_0\in \mathbb{N}$.
Thus, it suffices to show \eqref{t.transf.cm.1} for
continuous $\phi$ with $\Lambda(B_{\eta(\kappa_\phi)})<1$.

Assume that $\phi\in \mathcal{L}_2$ is continuous on
$[0,T]^2$ and $\Lambda(B_{\eta(\kappa_\phi)})<1$.
Then $\kappa_\phi$ is also continuous on $[0,T]^2$.
By Lemmas~\ref{l.kappa.phi.trace} and \ref{l.trace}, we 
obtain that 
\[
    \text{\rm tr} B_{\kappa_\phi}
    =\int_0^T [\text{\rm tr}\kappa_\phi(t,t)]dt
    =\int_0^T \biggl(\int_t^T
       [\text{\rm tr}\phi(t,s)]ds\biggr)dt
    =\int_0^T \biggl(\int_0^s[
         \text{\rm tr}\phi(t,s)]dt\biggr)ds.
\]
With the help of the identity
${\det}_2(I+A)=\det(I+A)e^{-\text{\rm tr} A}$ for
$A\in\mathcal{H}^{\otimes 2}$ of trace class,
plugged into \eqref{l.transf.cm.1}, this implies 
\eqref{t.transf.cm.1}.
\end{proof}

\begin{remark}\label{r.transf.cm}
We shall compare Theorem~\ref{t.transf.cm} with the change
of variables formula achieved by Cameron and Martin
\cite{CM}. 
To do so, we first recall their formula.
Let $d=1$ and $T=1$.
They considered the linear transformation of the form
\[
    \mathbb{T}:\mathcal{W}\ni w\mapsto
       w+\int_0^1 K(\cdot,s)w(s)ds,
\]
where $K\in \mathcal{L}_2$ and $K(0,\cdot)=0$. 
Reminding that their ``Wiener measure'' is the distribution
of the mapping $w\mapsto (1/\sqrt{2})w$ under $\mu$, we
rewrite their change of variables formula in our setting as
follows.  
Under several assumptions, which we do not restate here, it
holds that
\begin{equation}\label{eq:cm.1}
    |D|\int_S (f\circ \mathbb{T}) e^{-\Phi} d\mu
    =\int_{\mathbb{T}(S)} f d\mu
\end{equation}
for every Borel subset $S$ of $\mathcal{W}$ and 
$f\in C_b(\mathcal{W})$, where 
$\mathbb{T}(S)=\{\mathbb{T} w;w\in S\}$,
\begin{equation}\label{eq:cm.2}
    D=1+\sum_{n=1}^\infty \frac1{n!} \int_0^1\cdots\int_0^1
      \det\Bigl(\bigl(K(s_i,s_j)\bigr)_{1\le i,j\le n}\Bigr)
      ds_1\cdots ds_n,
\end{equation}
and
\begin{align*}
    \Phi(w)
    =& \frac12\int_0^1\biggl[\frac{d}{dt} 
          \int_0^1 K(t,s)w(s)ds\biggr]^2dt
     +\int_0^1\biggl[\int_0^1 
        \frac{\partial K}{\partial t}(t,s) w(s) ds
        \biggr]dw(t)
   \\
    & +\frac12\int_0^1 J(s) d\{w(s)^2\}
    \quad\text{with }
    J(s)=\lim_{t\downarrow s}K(t,s)
         -\lim_{t\uparrow s}K(t,s).
\end{align*}
They said that both ``$dw(t)$'' and ``$d\{w(s)^2\}$'' exist
as Stieltjes integrals.   

To compare \eqref{eq:cm.1} with \eqref{t.transf.cm.1},
we assume that $K$ is continuous on $[0,1]^2$ and of the
form 
\[
    K(t,s)=\int_0^t \phi(u,s)du
    \quad\text{for }(t,s)\in[0,1]^2
\]
for some $\phi\in\mathcal{L}_2$ with 
$\Lambda(B_{\kappa_\phi})<1$.
Then 
$\mathbb{T}=\iota+\mathbb{F}_\phi
           =\iota+F_{\kappa_\phi}$, and 
Theorem~\ref{t.transf.cm} is applicable to this $\phi$.
Combining the following observations (i)--(v), we see that
\eqref{eq:cm.1} follows from \eqref{t.transf.cm.1}.
\\
(i) The first term of $-\Phi$ is equal to the second one of 
$\Psi_\phi$.
\\
(ii)
Since $\partial K/\partial t=\phi$, exchanging the order of
the double integral with respect to $ds$ and ``$dw(t)$'',
we rewrite the second term of $-\Phi$ as 
\[
    -\int_0^1\biggl[\int_0^1\phi(t,s)dw(t)\biggr]
       w(s)ds.
\]
Since $d=1$, $M^\dagger=M$ for $M\in \mathbb{R}^{1\times 1}$
and the multiplications is commutative.
Hence this term coincides with the first term of
$\Psi_\phi$ by regarding ``$dw(t)$'' as the It\^o integral
$d\theta(t)$. 
\\
(iii) 
$J=0$ and the third term of $-\Phi$ vanishes.
\\
(iv)
For $\sigma\in\mathcal{L}_2$, define the linear operator
$\mathcal{L}_\sigma$ on $L^2([0,1];\mathbb{R})$ by 
\[
    \mathcal{L}_\sigma f=\int_0^1 \sigma(\cdot,s)f(s)ds
    \quad\text{for }f\in L^2([0,1];\mathbb{R}).
\]
We identify $\mathcal{L}_\sigma$ and $B_\sigma$ via the
correspondence between  $L^2([0,1];\mathbb{R})$ and
$\mathcal{H}$ such that 
$L^2([0,1];\mathbb{R})\ni f\leftrightarrow 
 \int_0^\bullet f(t)dt\in \mathcal{H}$.
Since 
\[
    K(t,s)=\int_0^1 \psi(t,u)\phi(u,s)du
    \quad\text{for }(t,s)\in[0,1]^2, 
\]
where 
$\psi(t,s)=\boldsymbol{1}_{[0,t)}(s)$ for 
$(t,s)\in[0,1]^2$ as before, 
$\mathcal{L}_K=\mathcal{L}_\psi \mathcal{L}_\phi$, 
and it is of trace class.
Due to the continuity of $K$, $\det(I+\mathcal{L}_K)$ admits
the Fredholm representation as described in
\eqref{eq:cm.2} (\cite[Theorem~VI.1.1]{GGK}).
Hence 
\[
    D=\det(I+\mathcal{L}_K)
     =\det(I+\mathcal{L}_\psi \mathcal{L}_\phi).
\]
If $A_1,A_2\in \mathcal{H}^{\otimes 2}$ are both of trace
class, then $\det(I+A_1A_2)=\det(I+A_2A_1)$ 
(\cite[IV.(5.9)]{GGK}).
This identity is extended to $\mathcal{H}^{\otimes 2}$ by a
standard approximation argument.
Thus we have that
\[
    \det(I+B_\psi B_\phi)=\det(I+B_\phi B_\psi).
\]
Hence, by \eqref{l.kappa.phi.trace.21} and the
identification  between $\mathcal{L}_\sigma$ and $B_\sigma$,
we obtain that 
\[
    D=\det(I+B_{\kappa_\phi}).
\]
(v)
Let $S$ be a Borel subset of $\mathcal{W}$ and 
$f\in C_b(\mathcal{W})$.
Put
$\Phi=\boldsymbol{1}_S\circ
    (\iota+F_{\widehat{\kappa_\phi}})$.
By Theorem~\ref{t.inv.transf}, it holds that 
\[
     f(\iota+F_{\kappa_\phi})
       \boldsymbol{1}_S
     =\bigl(f \Phi\bigr)(\iota+F_{\kappa_\phi})
    \quad\text{and}\quad
    f \Phi
    = f \boldsymbol{1}_{\mathbb{T}(S)}
    \quad\text{$\mu$-a.s.}
\]
\end{remark}

\begin{remark}\label{r.gen.cv}
We compare \eqref{eq:transf} with the identity obtained from
the general change of variables formula described in
\cite[Theorem~5.6.1]{mt-cambridge}.
The change of variables formula asserts that $F\in
\mathbb{D}^\infty(\mathcal{H})$ having an $r\in(1/2,\infty)$
with 
$\exp(-D^*F+r\|DF\|_{\mathcal{H}^{\otimes 2}}^2)
 \in L^{1+}(\mu)\equiv \bigcup_{p\in(1,\infty)}L^p(\mu)$ 
satisfies that 
\[
    \int_{\mathcal{W}} f(\iota+F)
      {\det}_2(I+DF)e^{-D^*F-(\|F\|_{\mathcal{H}}^2/2)} d\mu
    =\int_{\mathcal{W}} f d\mu
\]
for every $f\in C_b(\mathcal{W})$.
Let $\kappa\in \mathcal{L}_2$.
Applying this assertion to $F_\kappa$ with the help of
\eqref{l.f.kappa.21} and Lemma~\ref{l.d*f.kappa}, we obtain
that $F_\kappa$ with $\exp(-D^*F_\kappa)\in L^{1+}(\mu)$ 
satisfies that  
\[
    {\det}_2(I+B_\kappa)\int_{\mathcal{W}} f(\iota+F_\kappa)
      e^{\mathfrak{q}_{\eta(\kappa)}} d\mu
    =e^{\|\kappa\|_2^2/2}\int_{\mathcal{W}} f d\mu
\]
for every $f\in C_b(\mathcal{W})$.
Substituting $f=1$ into this, we see that 
${\det}_2(I+B_\kappa)>0$, and hence arrive at
\eqref{eq:transf}.
By \eqref{l.d*f.kappa.21}, we have that
\[
    -D^*F_\kappa=\mathfrak{q}_{s(\kappa)},
\]
where $s(\kappa)\in \mathcal{S}_2$ is given by
\[
    s(\kappa)(t,s)=-\{\kappa(t,s)+\kappa(s,t)^\dagger\}
    \quad\text{for }(t,s)\in[0,T]^2.
\]
Since $p \mathfrak{q}_{s(\kappa)}=\mathfrak{q}_{ps(\kappa)}$
and  $pB_{s(\kappa)}=B_{ps(\kappa)}$ for any
$p\in(1,\infty)$, 
by Lemma~\ref{l.q.eta.int}, we see that
$\exp(-D^*F_\kappa)\in L^{1+}(\mu)$
if and only if $\Lambda(B_{s(\kappa)})<1$. 
Thus the general change of variables formula implies 
\eqref{eq:transf} with ${\det}_2(I+B_\kappa)>0$ if 
$\Lambda(B_{s(\kappa)})<1$.

Assume that $\Lambda(B_{s(\kappa)})<2$.
Then $\Lambda(B_{\eta(\kappa)})<1$ and 
${\det}_2(I+B_\kappa)>0$.
In fact, it holds that
\[
    \inf_{\|h\|_{\mathcal{H}}=1}
      \langle B_{a\kappa} h,h\rangle_{\mathcal{H}}
    =-\frac{a}2 \Lambda(B_{s(\kappa)})>-a
    \quad\text{for }a\in[0,1].
\]
Hence 
\[
    \inf_{\|h\|_{\mathcal{H}}=1}
       \|(I+B_{a\kappa})h\|_{\mathcal{H}}
    \ge \inf_{\|h\|_{\mathcal{H}}=1}
        \langle (I+B_{a\kappa})h,h\rangle_{\mathcal{H}}
    >0
    \quad\text{for }a\in[0,1].
\]
By Remark~\ref{r.transf}(i), 
$\Lambda(B_{\eta(a\kappa)})<1$ for $a\in[0,1]$.
Due to the observation made before
Theorem~\ref{t.inv.transf}, this implies that
${\det}_2(I+B_{a\kappa})\ne0$ for $a\in[0,1]$.
Since the mapping 
$[0,1]\ni a\mapsto 
 {\det}_2(I+B_{a\kappa})={\det}_2(I+aB_\kappa)
 \in \mathbb{R}$ is
continuous and ${\det}_2(I+B_{a\kappa})=1$ for $a=0$, 
${\det}_2(I+B_{a\kappa})>0$ for $a\in[0,1]$. 
In particular, $\Lambda(B_{\eta(\kappa)})<1$ and 
${\det}_2(I+B_{\kappa})>0$.
Thus Theorem~\ref{t.transf} covers the cases obtained by 
using the general change of variables formula. 

In general, even if $\Lambda(B_{\eta(\kappa)})<1$ and 
${\det}_2(I+B_{\kappa})>0$, it does not hold that 
$\Lambda(B_{s(\kappa)})<2$.
For example, take orthonormal $h_1,h_2\in \mathcal{H}$ and 
$b_1,b_2<-1$.
Define $\kappa\in \mathcal{S}_2$ by 
\[
    \kappa(t,s)=b_1[h_1^\prime(s)\otimes h_1^\prime(t)]
     +b_2[h_2^\prime(s)\otimes h_2^\prime(t)]
    \quad\text{for }(t,s)\in[0,T]^2.
\]
We have that
\begin{align*}
    & s(\kappa)(t,s)
      =-2b_1[h_1^\prime(s)\otimes h_1^\prime(t)]
       -2b_2[h_2^\prime(s)\otimes h_2^\prime(t)],
    \\
    & 
    \eta(\kappa)(t,s)
       =-(2b_1+b_1^2)[h_1^\prime(s)\otimes h_1^\prime(t)]
        -(2b_2+b_2^2)[h_2^\prime(s)\otimes h_2^\prime(t)]
    \quad\text{for }(t,s)\in[0,T]^2.
\end{align*}
Hence 
\[
    \Lambda(B_{s(\kappa)})=-2(b_1\wedge b_2)>2
    \quad\text{and}\quad
    \Lambda(B_{\eta(\kappa)})=1-\{1+(b_1\vee b_2)\}^2
    <1.
\]
where $b_1\wedge b_2=\min\{b_1,b_2\}$.
Further, 
\[
    {\det}_2(I+B_\kappa)=(1+b_1)(1+b_2)e^{-(b_1+b_2)}>0.
\]

Thus, our condition that $\Lambda(B_{\eta(\kappa)})<1$
covers a wider class of transformations of order one than
the class obtained via the general change of variables
formula.
\end{remark}



\end{document}